\documentclass[11pt]{article}
\usepackage{pxfonts}
\usepackage{yfonts}
\usepackage{dsfont}
\usepackage{graphicx}
\usepackage{relsize}
\usepackage{color}

\def\bel{\begin{equation}\label}
\def\eeq{\end{equation}}
\def\ds{\displaystyle}

\def\mt{\longrightarrow}
\def\v{\vskip 1em}
\def\vsk{\vskip 40em}
\def\ve{\varepsilon}
\def\R{\mathbb R}
\def\Z{\mathbb Z}
\def\C{\mathfrak{B}}
\def\Cx{\mathds C}
\def\N{{\bf N}}
\def\exp{{\bf exp}}
\def\Re{{\bf Re}}
\def\Im {{\bf Im}}
\def\S{{\bf S}}

\def\F{{\bf F}}

\def\O{{\bf O}}
\def\Q{{\bf Q}}
\def\D{{\bf D}}
\def\J{{\bf J}}

\def\B{{\bf B}}
\def\H{{\bf H}}
\def\L{{\bf L}}
\def\T{{\bf T}}

\def\BMO{{\bf BMO}}

\def\p{{\partial}}
\def\a{{\bf a}}
\def\b{{\bf b}}
\def\i{{\bf i}}

\def\Hat{\widehat}

\def\bar{\overline}
\def\supp{{\bf supp}}

\def\I{{\bf I}}

\def\M{{\bf M}}

\def\vol{{\bf vol}}

\def\alpha{\alphaup}
\def\beta{\betaup}
\def\gamma{\gammaup}
\def\delta{\deltaup}
\def\theta{\thetaup}
\def\xi{{\xiup}}
\def\eta{{\etaup}}
\def\tau{{\tauup}}
\def\rho{{\rhoup}}
\def\phi{{\phiup}}
\def\psi{{\psiup}}
\def\lambda{{\lambdaup}}
\def\omega{\omegaup}
\def\varphi{{\varphiup}}
\def\gamma{{\gammaup}}
\def\c{{\bf c}}

\def\t{{\bf t}}

\def\m{{\bf m}}

\def\u{{\bf u}}

\newtheorem{remark}{Remark}[section]

\begin{document}
 \[\begin{array}{cc}\hbox{\LARGE{\bf A multi-parameter family of Fourier integral operators}}
  \end{array}\]
 
 \[\hbox{Mengmeng Dou~~~~~~~~Zipeng Wang~~~~~~~~Jiashu Zhang}\]
\begin{abstract}
We study a new class of Fourier integral operators defined in $\R^\N$. Their  symbols  are allowed to satisfy a differential inequality with certain multi-parameter characteristic. We prove these operators of order $-{\N-1\over 2}$   bounded from the classical, atom decomposable $\H^1$-Hardy space to $\L^1(\R^\N)$. As a result, we obtain a sharp $\L^p$-estimate. 
\v
Simultaneously, a generalized Sobolev $\L^p$-space is introduced. We establish the   Sobolev $\L^p$-norm inequality for convolutions with a distribution having singularity on the unit sphere. As an application, we give a new a priori estimate for  the solution of wave equations by requiring less regularity on  the source term and initial data. 
\end{abstract}

\section{Introduction}
\setcounter{equation}{0}
In this paper, we revisit on a classical problem analyzing the propagation of waves.   Denote $\Gamma$ and $\J$ for Gamma and Bessel functions. As first investigated by Strichartz \cite{Strichartz}, $\Omega^\delta$ is a distribution defined  in $\R^\N$ by analytic continuation from 
 \bel{Omega^delta}
\begin{array}{ccc}\ds
\Omega^\delta(x)~=~
\pi^{\delta-1}\Gamma^{-1}(\delta)\left({1\over 1-|x|^2}\right)^{1-\delta}_+,\qquad \hbox{\small{$\Re\delta>0$}}.
\end{array}
\eeq 
Alternatively, $\Omega^\delta$ can be defined by its Fourier transform 
\bel{Omega^delta Transform}
\begin{array}{lr}\ds
\Hat{\Omega}^\delta(\xi)~=~\left({1\over |\xi|}\right)^{{\N-1\over 2}+\delta-{1\over 2}}\J_{{\N-1\over 2}+\delta-{1\over 2}}\left(2\pi|\xi|\right),\qquad \hbox{\small{$\delta\in\Cx$}}.
\end{array}
\eeq
See appendix A in the end.

Sobolev  $\L^p$-norm inequalities  regarding $f\ast\Omega^\delta$  have been systematically studied, for example by Strichartz \cite{Strichartz}, Miyachi \cite{Miyachi},  Peral \cite{Peral} and Rubin \cite{Rubin}.
Essentially, it is equivalent to assert the $\L^p$-regularity of Fourier integral operators:
\bel{Ff}
\F f(x)~=~\int_{\R^\N} e^{2\pi\i \Phi(x,\xi)} \sigma(\xi) \Hat{f}(\xi)d\xi,\qquad \Phi(x,\xi)~=~x\cdot\xi\pm |\xi|
\eeq
where $\sigma\in\mathcal{C}^\infty(\R^\N)$ is vanished inside the unit ball.

$\diamond$ {\small Throughout, $\C>0$ is regarded as a generic constant depending on its sub-indices}.

We say $\sigma\in\S^{-\m}$ if 
\bel{class}
\left|\partial^\alpha_\xi \sigma(\xi)\right|~\leq~\C_\alpha~\left({1\over 1+|\xi|}\right)^{\m+|\alpha|}
\eeq
for every multi-index $\alpha$. 

Suppose $\sigma\in\S^0$. Plancherel theorem implies that $\F$  is bounded on $\L^2(\R^\N)$. On the other hand, $\F$ of order zero is not bounded on $\L^p(\R^\N)$ if $p\neq2$. Choose $\sigma(\xi)=(1+|\xi|)^{-\m}$. Then, $\F$ is not bounded on $\L^p(\R^\N)$ for $\left|{1\over p}-{1\over 2}\right|>{\m\over \N-1}$. See {\bf 6.13}, chapter IX of Stein \cite{Stein}. The following remarkable result is due to Miyachi \cite{Miyachi}.
\v
{\bf Theorem A:~~Miyachi, 1980}~~
{\it Let $\F$ defined in (\ref{Ff}). Suppose $\sigma\in \S^{-{\N-1\over 2}}$. We have
	\bel{Result A H}
	\left\| \F f\right\|_{\L^1(\R^\N)} ~\leq~\C~\left\| f\right\|_{\H^1(\R^\N)}.
	\eeq 
	Suppose $\sigma\in \S^{-\m}$ for $0\leq\m<(\N-1)/2$. We have
	\bel{Result A}
	\begin{array}{cc}\ds
		\left\| \F f\right\|_{\L^p(\R^\N)}~\leq~\C_p~\left\| f\right\|_{\L^p(\R^\N)},\qquad 
		\left|{1\over p}- {1\over 2}\right|~\leq~{\m\over \N-1}.
	\end{array}
	\eeq}
\begin{remark}  $\H^1(\R^\N)$ in (\ref{Result A H}) is the  $\H^1$-Hardy space  investigated by Fefferman and Stein \cite{Fefferman-Stein}.
\end{remark}

Denote $\L^p_r(\R^\N)$ to be the Sobolev $\L^p$-space of order $r\in\R$ for $1\leq p< \infty$. 
By using the $\L^p$-estimate in (\ref{Result A}), Rubin \cite{Rubin} proved the next Sobolev $\L^p$-norm inequality.

\v
{\bf Theorem B:~~Rubin, 1991}~~{\it Let $\Omega^\delta$  defined as (\ref{Omega^delta})-(\ref{Omega^delta Transform}). 
	For $-{\N-1\over2}\leq\delta-r+s<0$, we have  
	\bel{Result B}
	\left\|  f\ast\Omega^\delta\right\|_{\L^p_r(\R^\N)}~\leq~\C_{p~\delta~r~s}~\left\| f\right\|_{\L^p_s(\R^\N)},\qquad
	\left|{1\over p}-{1\over 2}\right|~\leq~{\delta-r+s\over \N-1}+{1\over 2}.
	\eeq}
\begin{remark} 
 The restriction on the range of $p$ inside (\ref{Result B}) is also  necessary.  Moreover, a  complete characterization of $f\ast\Omega^\delta\colon\L^p_s(\R^\N)\mt\L^q_r(\R^\N)$ for $\delta-r+s\in\R$ and $1\leq p\leq q\leq\infty$ has been established. See p.403 in the paper of Rubin \cite{Rubin}. 
\end{remark}
A direct application of {\bf Theorem B} goes to the Cauchy problem for the wave equation
\bel{Wave}
\begin{array}{cc}\ds
\partial^2_t \u(x,t)-\Delta_x \u(x,t)~=~h(x,t),\qquad (x,t)\in\R^\N\times(0,T],
\\\\ \ds
  \u(x,0)~=~f(x),\qquad  \partial_t \u(x,0)~=~g(x).
\end{array}
\eeq
As a well known result, we have
\[
\u(x,t)~=~f\ast \Omega^{-{\N-1\over2}}_t(x)+t g\ast\Omega^{-{\N-3\over 2}}_t(x)+\int_0^t h\ast\Omega^{-{\N-3\over2}}_\tau (x,\tau) (t-\tau)d\tau
\]
where $\Hat{\Omega}^{-{\N-1\over 2}}(\xi)=\cos(2\pi|\xi|)$ and $\Hat{\Omega}^{-{\N-3\over 2}}(\xi)={\sin(2\pi|\xi|)\over |\xi|}$.

The Sobolev $\L^p$-norm inequality in (\ref{Result B}) implies
\bel{u norm}
\begin{array}{cc}\ds
\left\| \u(\cdot, t)\right\|_{\L^p_r(\R^\N)}~\leq~\C_{t~p~r~s} ~\left\{\left\| f\right\|_{\L^p_s(\R^\N)}+\left\| g\right\|_{\L^p_{s-1}(\R^\N)}+\int_0^t \left\| h(\cdot, t-\tau)\right\|_{\L^p_{s-1}(\R^\N)} d\tau\right\}
\\\\ \ds
\hbox{whenever}\qquad  \left|{1\over p}-{1\over 2}\right|~\leq~{s-r\over \N-1}.
\end{array}
\eeq
Recall $\F$ defined in (\ref{Ff}).
To handle more  hyperbolic equations consisting of variable coefficients, a generalization occurred.
The phase 
$\Phi(x,\xi)=x\cdot\xi\pm |\xi|$ is replaced by any real-valued  function which is homogeneous of degree $1$ in $\xi$ and smooth in $(x,\xi)$ for $\xi\neq0$. Moreover, it satisfies   
 the non-degeneracy condition:
\[\det\left[{\p^2\Phi\over \p x\p \xi}\right]\left(x,\xi\right)\neq0, \qquad\xi\neq0.\] 
On the other hand, the symbol becomes $\sigma(x,\xi)\in\mathcal{C}^\infty(\R^\N\times\R^\N)$ having a $x$-compact support. The optimal $\L^p$-estimate regarding  Fourier integral operators defined within this general setting was first investigated by H\"{o}rmander \cite{Hormander},  
Duistermaat and H\"{o}rmander  \cite{Duistermaat-Hormander} and then by  Colin de Verdi\'{e}re and Frisch \cite{Colin-Frisch},  Brenner \cite{Brenner}, Peral \cite{Peral}, Miyachi \cite{Miyachi}, Beals \cite{Beals} and eventually proved by Seeger, Sogge and Stein \cite{Seeger-Sogge-Stein}.

Over the several past decades, the theory of Fourier integral operators has been extensively developed. 
For instance, a degenerate Fourier integral operator  is introduced by Seeger \cite{Seeger}. Certain bilinear Fourier integrals have invented by Rodriguez-L\'{o}pez, Rule and Staubach \cite{Lopez-Rule-Staubach} and 
Kato, Miyachi and Tomita \cite{Kato-Miyachi-Tomita}. A class of bi-parameter Fourier integral operators has investigated by Hong, Lu and Zhang \cite{Hong-Lu-Zhang}. Furthermore, see Ruzhansky and Sugimoto \cite{Ruzhansky-Sugimoto} and Sindayigaya \cite{Sindayigaya} for  the regarding global estimates. More recently, Wu and Yang \cite{Wu-Yang} proves the $\L^p$-boundedness for rough Fourier integral operators.

Next, we give a multi-parameter extension  to both {\bf Theorem A}  and {\bf Theorem B}.

\subsection{Multi-parameter extensions}
The study of certain operators  that  commute with a  family of multi-parameter dilations,  dates back to the time of  Jessen, Marcinkiewicz and Zygmund.  
The $\L^p$-theory for the strong maximal function and singular integrals defined on product spaces have been established. For instance,  see the references by Cordoba and Fefferman \cite{Cordoba-Fefferman}, Fefferman and Stein \cite{R.Fefferman-Stein}, Journ\'{e} \cite{Journe'},  Fefferman \cite{R.Fefferman}   and 
M\"{u}ller, Ricci and Stein \cite{Muller-Ricci-Stein}. 

In particular, for $0<p\leq1$, a bi-parameter  Hardy space, denoted by $\H^p\times \H^p(\R^{\N_1}\times\R^{\N_2})$ was introduced by  M.~Malliavin and P.~Malliavin \cite{Malliavin} and Gundy and Stein \cite{Gundy-Stein}. Such a Hardy space defined by the corresponding bi-parameter Littlewood-Paley inequality  cannot be characterized in terms of  "rectangle atoms". See the counter-example of Carleson \cite{Carleson}. Surprisingly, by testing on any single rectangle atom and using a geometric covering lemma due to Journ\'{e} \cite{Journe'}, Fefferman \cite{R.Fefferman} is able to conclude that a bi-parameter Calder\'{o}n-Zygmund operator defined on $\R^{\N_1}\times\R^{\N_2}$ is bounded from $\H^p\times \H^p(\R^{\N_1}\times\R^{\N_2})$ to $\L^p(\R^{\N_1}\times\R^{\N_2})$. 

Consider $\N=\N_1+\N_2+\cdots+\N_n$ where $n\ge2$ is the number of parameters. We write $\xi=(\xi_1,\xi_2,\ldots,\xi_n)^T\in\R^{\N_1}\times\R^{\N_2}\times\cdots\times\R^{\N_n}=\R^\N$. Denote $\rho=(\rho_1,\rho_2,\ldots,\rho_n)$ for which $\rho_i\ge0, i=1,2\ldots,n$. 
We say $\sigma\in\S^{-\m}_\rho$ if
\bel{Class}
\begin{array}{cc}\ds
\left|\p_\xi^\alpha \sigma(\xi)\right|~\leq~\C_\alpha~ \prod_{i=1}^n \left({1\over 1+|\xi_i|}\right)^{|\alpha_i|+\rho_i} 
\\\\ \ds
\m=\rho_1+\rho_2+\cdots+\rho_n,\qquad \rho_i=\m\theta_i,\qquad\hbox{${\N_i-1\over \N-1}<\theta_i<{\N_i\over \N-1}$},\qquad  i=1,2,\ldots,n.
\end{array}
\eeq
for every multi-index $\alpha$.
\v
{\bf Theorem One}~~
{\it Let $\F$ defined in (\ref{Ff}). Suppose $\sigma\in \S^{-{\N-1\over 2}}_\rho$. We have
	\bel{Result One H}
	\left\| \F f\right\|_{\L^1(\R^\N)} ~\leq~\C_\rho~\left\| f\right\|_{\H^1(\R^\N)}.
	\eeq 
	Suppose $\sigma\in \S^{-\m}_\rho$ for $0\leq\m<(\N-1)/2$. We have
	\bel{Result One}
	\begin{array}{cc}\ds
		\left\| \F f\right\|_{\L^p(\R^\N)}~\leq~\C_{\rho~p}~\left\| f\right\|_{\L^p(\R^\N)},\qquad 
		\left|{1\over p}- {1\over 2}\right|~\leq~{\m\over \N-1}.
	\end{array}
	\eeq}
\begin{remark}  $\H^1(\R^\N)$ in (\ref{Result One H}) is the same $\H^1$-Hardy space  investigated by Fefferman and Stein \cite{Fefferman-Stein}.
	Furthermore, it has a characterization of atom decomposition established by Coifman \cite{Coifman}.
\end{remark}
In contrast to the $\H^p\times \H^p\mt\L^p$-estimate for the bi-parameter Calder\'{o}n-Zygmund operator, 
$\F$ defined in (\ref{Ff}) with $\sigma\in\S^{-{\N-1\over 2}}_\rho$ is an example of multi-parameter operators  bounded from the classical, atom decomposable $\H^1(\R^\N)$ to $\L^1(\R^\N)$.

Now, we introduce a generalized  Sobolev $\L^p$-space having a multi-parameter characteristic. 
Let $\B_{\rho~s}$ be a distribution whose Fourier transform equals
\bel{B Transform}
\begin{array}{cc}\ds
\Hat{\B}_{\rho~s}(\xi)~=~\Big(1+|\xi|^2\Big)^{s\over 2}\prod_{i=1}^n \left(1+|\xi_i|^2\right)^{\rho_i\over2},\qquad s\in\R,
\\\\ \ds
\m=\rho_1+\rho_2+\cdots+\rho_n,\qquad \rho_i=\m\theta_i,\qquad\hbox{${\N_i-1\over \N-1}<\theta_i<{\N_i\over \N-1}$},\qquad  i=1,2,\ldots,n.
\end{array}
\eeq
Define
\bel{Sobolev space}
\L^p_{\rho~s}\left(\R^\N\right)~=~\Bigg\{f\colon f\ast\B_{\rho~s}\in\L^p\left(\R^\N\right)\Bigg\},\qquad 1< p<\infty.
\eeq
\begin{remark}  $\L^p_{\rho~s}(\R^\N)\supset\L^p_{\m+s}(\R^\N)$. 
\end{remark}
This can be seen by considering $\ds\left(1+|\xi|^2\right)^{-{\m\over2}}\prod_{i=1}^n \left(1+|\xi_i|^2\right)^{\rho_i\over 2}$, $\m=\rho_1+\rho_2+\cdots+\rho_n$ which is a $\L^p$-Fourier multiplier for $1<p<\infty$. 
\begin{remark}
$\L^p_{\m}(\R^{\N})\subsetneq \L^p_{\rho~0}(\R^{\N})$ for $\m=\rho_1+\rho_2+\cdots+\rho_n>0$.
\end{remark}
Choose $f(x)=f_1(x_1)f_2(x_2)\cdots f_n(x_n)$ of which $f_i\in \L^p_{\rho_i}(\R^{\N_i})$: $\left[\left(1+|\xi_i|^2\right)^{\rho_i\over 2}\Hat{f}_i(\xi_i)\right]^\vee\in\L^p(\R^{\N_i})$ but 
$f_i\notin \L^p_\m(\R^{\N_i})$ as $\rho_i<\m$ for every $i=1,2,\ldots,n$. We find $f\in\L^p_{\rho~0}(\R^\N)$ and $f\notin\L^p_{\m}(\R^\N)$.
\v

{\bf Theorem Two} ~~
{\it  Let $\Omega^\delta$  defined as (\ref{Omega^delta})-(\ref{Omega^delta Transform}). For $-{\N-1\over2}\leq\delta-r+s+\m<0$, we have 
\bel{Result Two}
\left\| f\ast\Omega^\delta\right\|_{\L^p_r(\R^\N)}~\leq~\C_{p~\delta~r~s~\rho}~\left\| f\right\|_{\L^p_{\rho~s}(\R^\N)}
,\qquad \left|{1\over p}-{1\over 2}\right|~\leq~{\delta-r+s+\m\over \N-1}+{1\over 2}.
\eeq}

Recall the inhomogeneous wave equation in (\ref{Wave}). Moreover, its solution is given as
\[
\u(x,t)~=~f\ast \Omega^{-{\N-1\over2}}_t(x)+t g\ast\Omega^{-{\N-3\over 2}}_t(x)+\int_0^t h\ast\Omega^{-{\N-3\over2}}_\tau (x,\tau) (t-\tau)d\tau.
\]
For $-{\N-1\over2}\leq\delta-r+s+\m<0$, the Sobolev $\L^p$-norm inequality in (\ref{Result Two}) implies
\bel{u norm rho}
\begin{array}{cc}\ds
\left\| \u(\cdot, t)\right\|_{\L^p_r(\R^\N)}~\leq~\C_{t~p~r~s~\rho} ~\left\{\left\| f\right\|_{\L^p_{\rho~s}(\R^\N)}+\left\| g\right\|_{\L^p_{\rho~s-1}(\R^\N)}+\int_0^t \left\| h(\cdot, t-\tau)\right\|_{\L^p_{\rho~s-1}(\R^\N)} d\tau\right\}
\\\\ \ds
\hbox{whenever}\qquad  \left|{1\over p}-{1\over 2}\right|~\leq~{\m+s-r\over \N-1}.
\end{array}
\eeq
Compare (\ref{u norm rho}) to (\ref{u norm}). Let $f\in\L^p_\m(\R^\N)$, $g\in \L^p_{\m-1}(\R^\N)$ and $h(\cdot, t)\in \L^p_{\m-1}(\R^\N), 0< t\leq T$ where $\m=\rho_1+\rho_2+\cdots\rho_n>0$. From (\ref{u norm}), we find
\bel{u norm 0}
\begin{array}{cc}\ds
\left\| \u(\cdot, t)\right\|_{\L^p(\R^\N)}~\leq~\C_{t~p~\m} ~\left\{\left\| f\right\|_{\L^p_\m(\R^\N)}+\left\| g\right\|_{\L^p_{\m-1}(\R^\N)}+\int_0^t \left\| h(\cdot, t-\tau)\right\|_{\L^p_{\m-1}(\R^\N)} d\tau\right\}
\\\\ \ds
\hbox{whenever}\qquad  \left|{1\over p}-{1\over 2}\right|~\leq~{\m\over \N-1}.
\end{array}
\eeq
On the other hand, (\ref{u norm rho}) shows
\bel{u norm rho 0}
\begin{array}{cc}\ds
\left\| \u(\cdot, t)\right\|_{\L^p(\R^\N)}~\leq~\C_{t~p~\rho} ~\left\{\left\| f\right\|_{\L^p_{\rho~0}(\R^\N)}+\left\| g\right\|_{\L^p_{\rho~-1}(\R^\N)}+\int_0^t \left\| h(\cdot, t-\tau)\right\|_{\L^p_{\rho~-1}(\R^\N)} d\tau\right\}
\\\\ \ds
\hbox{whenever}\qquad  \left|{1\over p}-{1\over 2}\right|~\leq~{\m\over \N-1}.
\end{array}
\eeq
Recall {\bf Remark 1.5}. Observe that the same $\L^p$-norm of $\u(\cdot, t)$ can be controlled  by requiring less differentiability properties on $f,g$ and $h$.

In the next section, we show {\bf Theorem One} implying  {\bf Theorem Two}. Section 3,4 and 5 are devoted to the proof of {\bf Theorem One}
\subsection{Sketch on the proof of Theorem One}

In section 3, we develop a new framework where the frequency space is decomposed into an infinitely many dyadic cones. Each consisting partial operator has a Fourier transform supported on  these dyadic cones. Moreover, its kernel  enjoys certain majorization properties, accumulated as {\bf Lemma One}.  We prove {\bf Theorem One} by applying {\bf Lemma One}.

In section 4, we introduce a second dyadic decomposition. The ingenious construction was  initially given by Fefferman \cite{Fefferman'} and later refined by Seeger, Sogge and Stein \cite{Seeger-Sogge-Stein} in their study of Fourier integral operators. 

In order to prove {\bf Lemma One}, we need to combine the two different frameworks. The proof of {\bf Lemma One} is completed in section 5.

\vsk

\section{Proof of Theorem Two}
\setcounter{equation}{0}
Recall $\Hat{\Omega}^\delta(\xi)$ given in (\ref{Omega^delta Transform}) and
$\Hat{\B}_{\rho~s}(\xi)$ defined in (\ref{B Transform}). Let $-{\N-1\over 2}\leq\delta-r+s+\m<0$ where $\m=\rho_1+\rho_2+\cdots\rho_n$ and $\rho_i=\m\theta_i$ for ${\N_i-1\over \N-1}<\theta_i<{\N_i\over \N-1}$, $i=1,2,\ldots,n$. We have
\bel{Omega r}
\begin{array}{lr}\ds
\left(1+|\xi|^2\right)^{r\over2} \Hat{\Omega}^\delta(\xi) \left[\Hat{\B}_{\rho~s}(\xi)\right]^{-1}
~=~
\\\\ \ds
\left(1+|\xi|^2\right)^{r-s\over2}\left({1\over |\xi|}\right)^{{\N-1\over 2}+\delta-{1\over 2}}\J_{{\N-1\over 2}+\delta-{1\over 2}}\left(2\pi|\xi|\right) \prod_{i=1}^n \left({1\over 1+|\xi_i|^2}\right)^{\rho_i\over2}.
\end{array}
\eeq
Suppose $-{\N-1\over 2}-\kappa<\delta\leq -{\N-1\over 2}-\kappa+1$ for $\kappa=0,1,2\ldots$. By using the recurrence formula of Bessel functions in (\ref{J identity}), we write
\bel{J Sum}
\begin{array}{lr}\ds
\J_{{\N-1\over 2}+\delta-{1\over 2}}\left(2\pi|\xi|\right)~=~\pi^{-\kappa} \hbox{$\Big[ {\N-1\over 2}+\delta-{1\over 2}+1\Big]\cdots \Big[ {\N-1\over 2}+\delta-{1\over 2}+\kappa\Big]$} |\xi|^{-\kappa} \J_{{\N-1\over 2}+\delta-{1\over 2}+\kappa}\left(2\pi|\xi|\right)
\\\\ \ds ~~~~~~~~~~~~~~~~~~~~~~~~~~~
~+~\hbox{\bf intermediate terms}~+~(-1)^\kappa \J_{{\N-1\over 2}+\delta-{1\over 2}+2\kappa}\left(2\pi|\xi|\right).
\end{array}
\eeq
Each {\bf intermediate term} equals 
\bel{intermediate term}
|\xi|^{-\ell} \J_{{\N-1\over 2}+\delta-{1\over 2}+m}\left(2\pi|\xi|\right),\qquad 0<\ell<\kappa,\qquad \kappa<m<2\kappa
\eeq
multiplied by some constant depending on $\delta$ and $\kappa$. 

Observe that every Bessel function in the right hand side of (\ref{J Sum}) has a sub-index strictly larger than $-{1\over 2}$. The reason to write such expansion is to allow us applying the formulae in (\ref{Bessel}) and  (\ref{J asymptotic})-(\ref{J O}).

Given $0\leq\ell\leq \kappa$ and $\kappa\leq m\leq 2\kappa$, we consider
\bel{principal term}
\left(1+|\xi|^2\right)^{r-s\over2}\left({1\over |\xi|}\right)^{{\N-1\over 2}+\delta-{1\over 2}}|\xi|^{-\ell}\J_{{\N-1\over 2}+\delta-{1\over 2}+m}\left(2\pi|\xi|\right) \prod_{i=1}^n \left({1\over 1+|\xi_i|^2}\right)^{\rho_i\over2}.
\eeq
Let $\varphi\in\mathcal{C}^\infty_o(\R)$ be a smooth bump function  such that $\varphi(t)=1$ if $|t|\leq1$ and $\varphi(t)=0$ for $|t|>2$.
Define
\bel{phi_j}
\phi(\xi)~=~\varphi(|\xi|),\qquad
 \phi_j(\xi)~=~\varphi\left[ 2^{-j} |\xi|\right]-\varphi\left[2^{-j+1}|\xi|\right],\qquad j\in\Z.
\eeq
By using the integral formula  in (\ref{Bessel}), we find
\bel{omega_o}
\begin{array}{lr}\ds
\omega_o(\xi)~\doteq~\phi(\xi)\left(1+|\xi|^2\right)^{r-s\over2}|\xi|^{m-\ell}\left({1\over |\xi|}\right)^{{\N-1\over 2}+\delta-{1\over 2}+m}\J_{{\N-1\over 2}+\delta-{1\over 2}+m}\left(2\pi|\xi|\right) \prod_{i=1}^n \left({1\over 1+|\xi_i|^2}\right)^{\rho_i\over2}
\\\\ \ds~~~~~~~
~=~\hbox{$\pi^{{\N-1\over 2}+\delta-1+m}\Gamma^{-1}\left({\N-1\over 2}+\delta+m\right)$}
\\ \ds~~~~~~~~~~~~~~
\phi(\xi)|\xi|^{m-\ell} \left(1+|\xi|^2\right)^{r-s\over2}\prod_{i=1}^n \left({1\over 1+|\xi_i|^2}\right)^{\rho_i\over2}\int_{-1}^1 e^{2\pi\i|\xi| t } (1-t^2)^{{\N-1\over 2}+\delta-1+m} dt.
\end{array}
\eeq
This is a $\L^p$-Fourier multiplier for $1<p<\infty$ because it satisfies the Marcinkiewicz condition: 
$\left|\left(\xi\p_\xi\right)^\alpha \omega_o(\xi)\right|\leq\C_\alpha$
for every multi-index $\alpha$.

 By using the asymptotic expansion of Bessel functions in (\ref{J asymptotic})-(\ref{J O}), we write 
\bel{I rewrite}
\begin{array}{lr}\ds
\int_{\R^\N}e^{2\pi\i x\cdot\xi}\big[1-\phi(\xi)\big] \left(1+|\xi|^2\right)^{r-s\over2}\left({1\over |\xi|}\right)^{{\N-1\over 2}+\delta-{1\over 2}}|\xi|^{-\ell}\J_{{\N-1\over 2}+\delta-{1\over 2}+m}\left(2\pi|\xi|\right) \prod_{i=1}^n \left({1\over 1+|\xi_i|^2}\right)^{\rho_i\over2}\Hat{f}(\xi)d\xi
\\\\ \ds
~=~\I_{\ell~m}f(x)+\mathcal{E}_{\ell~m~N} f(x);
\end{array}
\eeq
\bel{II}
\begin{array}{lr}\ds
\I_{\ell~m}f(x)~=~\pi^{-1}\int_{\R^\N} e^{2\pi\i x\cdot\xi} \cos\left[2\pi|\xi|-\left({\N-1\over 2}+\delta-{1\over 2}+m\right) {\pi\over 2}-{\pi\over 4}\right]
\\ \ds~~~~~~~~~~~~~~~~~~~~~~~~~~~~~
\left(1+|\xi|^2\right)^{r-s\over2}|\xi|^{-\ell}\left({1\over |\xi|}\right)^{{\N-1\over 2}+\delta}\prod_{i=1}^n \left({1\over 1+|\xi_i|^2}\right)^{\rho_i\over2}[1-\phi(\xi)]
 \Hat{f}(\xi) d\xi
 \\\\ \ds
 ~+~\pi^{-1}  \sum_{k=1}^N \int_{\R^\N} e^{2\pi\i x\cdot\xi} \Bigg\{~  \a_k\left({1\over2\pi|\xi|}\right)^{2k} \cos\left[2\pi|\xi|-\left({\N-1\over 2}+\delta-{1\over 2}+m\right) {\pi\over 2}-{\pi\over 4}\right]
 \\ \ds~~~~~~~~~~~~~~~~~~~~~~~~~~~~~
+\b_k\left({1\over 2\pi|\xi|}\right)^{2k-1} \sin\left[2\pi|\xi|-\left({\N-1\over 2}+\delta-{1\over 2}+m\right) {\pi\over 2}-{\pi\over 4}\right]~ \Bigg\}
 \\ \ds~~~~~~~~~~~~~~~~~~~~~
\left(1+|\xi|^2\right)^{r-s\over2}|\xi|^{-\ell}\left({1\over |\xi|}\right)^{{\N-1\over 2}+\delta}\prod_{i=1}^n \left({1\over 1+|\xi_i|^2}\right)^{\rho_i\over2}[1-\phi(\xi)]
 \Hat{f}(\xi) d\xi
 \end{array}
 \eeq
where $\a_k,\b_k$ are constants given in (\ref{J asymptotic}).
 \bel{E}
  \mathcal{E}_{\ell~m~N} f(x)~=~ \int_{\R^\N} e^{2\pi\i x\cdot\xi} \hbox{\bf R}_N(\xi)  \left(1+|\xi|^2\right)^{r-s\over2}|\xi|^{-\ell}\left({1\over |\xi|}\right)^{{\N-1\over 2}+\delta}\prod_{i=1}^n \left({1\over 1+|\xi_i|^2}\right)^{\rho_i\over2}[1-\phi(\xi)]
 \Hat{f}(\xi) d\xi
 \eeq
 of which 
\[
\partial_\xi^\alpha \hbox{\bf R}_N(\xi)~=~\O\left(|\xi|^{-2N-{1\over 2}}\right),\qquad|\xi|\mt\infty
\]
for every multi-index $\alpha$.

The kernel of $\mathcal{E}_{\ell~m~N}$ belongs to $\L^1(\R^\N)$ for $N$   sufficiently large.
On the other hand, $\I_{\ell~m}$ defined in (\ref{II}) can be written as a finite sum of Fourier integral operators in (\ref{Ff}) with 
$\sigma(\xi)=\omega(\xi) \left[1-\phi(\xi)\right]$ for some $\omega\in\mathcal{C}^\infty(\R^\N\setminus 0)$ 
satisfying
\bel{Diff Ineq sigma}
\begin{array}{lr}\ds
\left|\partial_\xi^\alpha \sigma(\xi)\right|~\leq~\C_\alpha~ |\xi|^{-\ell}\left({1\over 1+ |\xi|}\right)^{{\N-1\over 2}+\delta-r+s} \prod_{i=1}^n \left({1\over 1+|\xi_i|}\right)^{\rho_i+|\alpha_i|}
\\\\ \ds~~~~~~~~~~~~~~
~\leq~\C_\alpha  \prod_{i=1}^n \left({1\over 1+|\xi_i|}\right)^{\left[{\N-1\over 2}+\delta-r+s\right]\theta_i+\rho_i+|\alpha_i|}\qquad \hbox{\small{( $|\xi|>1$, $\ell\ge0$ )}}
\\\\ \ds~~~~~~~~~~~~~~
~=~\C_\alpha  \prod_{i=1}^n \left({1\over 1+|\xi_i|}\right)^{\left[{\N-1\over 2}+\delta-r+s+\m\right]\theta_i+|\alpha_i|}\qquad\hbox{\small{( $\rho_i=\m\theta_i, i=1,2,\ldots,n$ )}}
\end{array}
\eeq
for every multi-index $\alpha$. 

Note that   $-{\N-1\over 2}\leq\delta-r+s+\m<0$.
By applying the $\L^p$-estimate in (\ref{Result One}), we find $\I_{\ell~m}f$ satisfying the Sobolev $\L^p$-norm inequality in (\ref{Result Two}).

The above estimates hold for every $0\leq\ell \leq \kappa$ and $\kappa\leq m\leq 2\kappa$. All together, we conclude {\bf Theorem Two}.

\section{Proof of Theorem One}
\setcounter{equation}{0}
Let $\varphi\in\mathcal{C}^\infty_o(\R)$ be the smooth bump function as before.
Denote $\t=(t_2,\ldots,t_n)$ for which $t_i, i=2,\ldots,n$ are non-negative integers.
We define
\bel{delta_t}
\begin{array}{cc}\ds
\delta_\t(\xi)~=~\prod_{i=2}^n \delta_{t_i}(\xi),
\qquad
 \delta_{t_i}(\xi)~=~ \varphi\left[2^{-t_i}{|\xi_1|\over|\xi_i|}\right]-\varphi\left[2^{-t_i+1}{|\xi_1|\over|\xi_i|}\right]
 \end{array}
\eeq
whose support is contained in the dyadic cone
\bel{Cone}
\hbox{V}_\t~=~ \Bigg\{  \xi\in\R^\N~\colon~      2^{t_i-1}~<~ {|\xi_1|\over|\xi_i|}~<~2^{t_i+1},~i=2,\ldots,n \Bigg\}.
\eeq

\begin{figure}[h]
\centering
\includegraphics[scale=0.40]{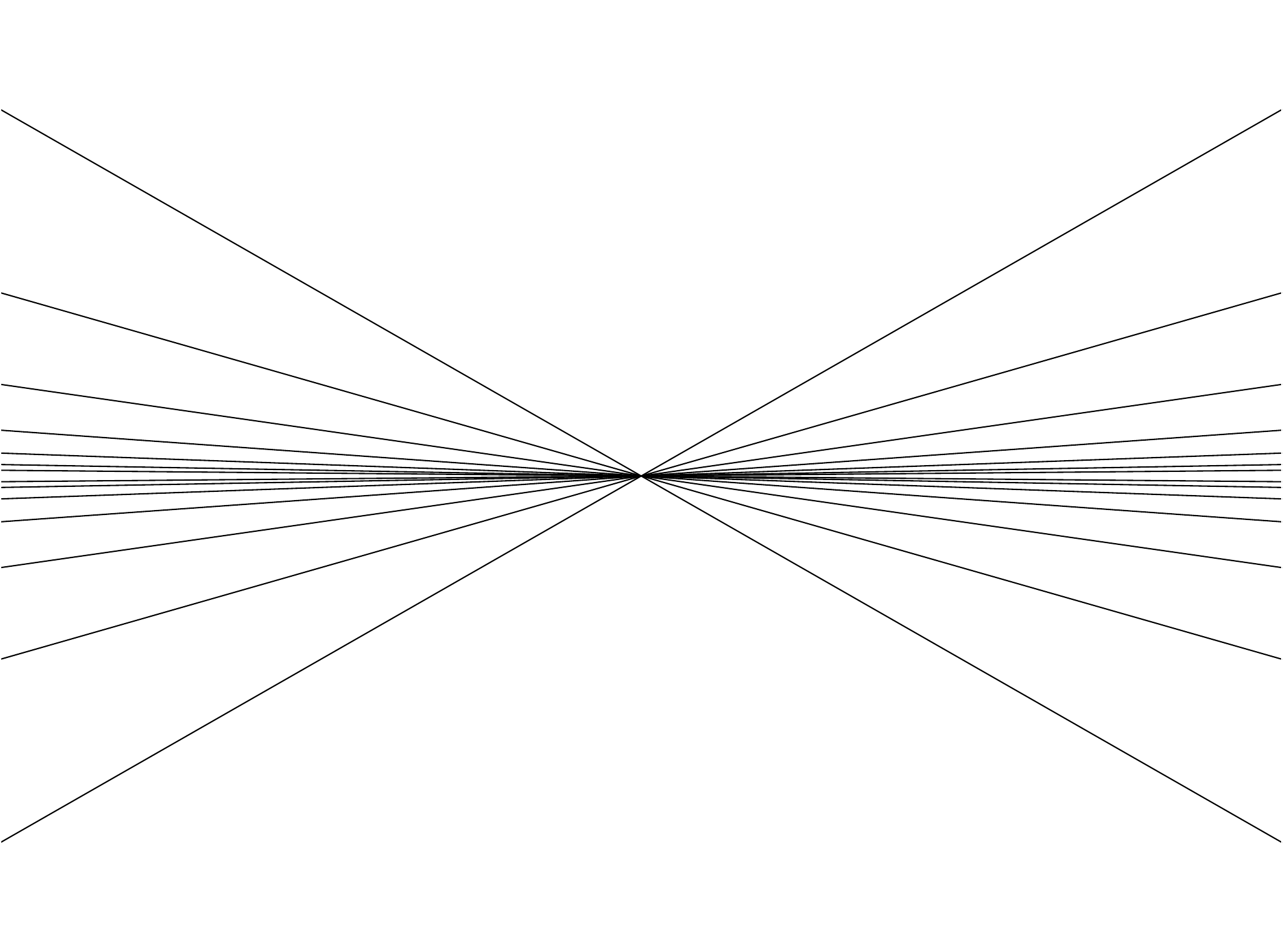}
\end{figure}

A  direct computation shows
\bel{delta_l Diff Ineq}
\left| \p_\xi^\alpha   \delta_\t(\xi)\right|~\leq~\C_\alpha~  |\xi_1|^{-|\alpha_1|}\prod_{i=2}^n |\xi_i|^{-|\alpha_i|}
\eeq
for every $(n-1)$-tuple $\t$ and every multi-index $\alpha$.

Consider
\bel{Delta}
\Delta(\xi)~=~\sum_{t_i\ge0,~ i=2,\ldots,n} \delta_\t(\xi)
\eeq
which has a support inside the conical region: $|\xi_1|\ge|\xi_i|, i=2,\ldots,n$.

Define
\bel{F tri}
\begin{array}{lr}\ds
\F^\triangle f(x)~=~\int_{\R^\N} e^{2\pi\i \Phi(x,\xi)} \Delta(\xi) \sigma(\xi)\Hat{f}(\xi) d\xi,\qquad \Phi(x,\xi)=x\cdot\xi\pm|\xi|.
\end{array}
\eeq
For symmetry reason, it is suffice to prove (\ref{Result One H})-(\ref{Result One}) for $\F^\triangle$ instead. 

Given $j>0$, we assert $I\cup II\cup III=\{2,\ldots,n\}$ such that
\bel{I,II}
0\leq t_i<j/2,~~~~ i\in I,\qquad  j/2\leq t_i<j,~~~~ i\in II,\qquad t_i\ge j,~~~~ i\in III.
\eeq
\begin{remark}
There are  $3^{n-1}$ such partitions: $I\cup II \cup III=\{2,\ldots,n\}$. 
\end{remark}
Let  $I\cup II\cup III$ to be fixed. Define
\bel{Delta_t j}
\Delta_{\t j~\left[I, \flat, \sharp\right]}(\xi)~=~ \mathfrak{S}^\flat_j(\xi)\mathfrak{S}^\sharp_j(\xi)\prod_{i\in I} \delta_{t_i}(\xi)
\eeq
where
\bel{S sharp flat}
\begin{array}{cc}\ds
\mathfrak{S}^\flat_j(\xi)~=~\prod_{i\in II}~\sum_{j/2\leq t_i< j} \delta_{t_i}(\xi),
\\\\ \ds 
\mathfrak{S}^\sharp_j(\xi) ~=~\prod_{i\in III}~\sum_{t_i\ge j} \delta_{t_i}(\xi)
~=~\prod_{i\in III} \left\{1-\varphi\left[2^{-j+1}{|\xi_1|\over|\xi_i|}\right]\right\}.
\end{array}
\eeq
\begin{remark}
Observe that
$\p_{\xi_1}^{\alpha_1}\p_{\xi_i}^{\alpha_i} \left[1-\varphi\left[2^{-j+1}{|\xi_1|\over|\xi_i|}\right]\right]=0$ whenever ${|\xi_1|\over |\xi_i|}<2^{j-1}$ or ${|\xi_1|\over |\xi_i|}>2^j$ and $|\alpha_1|+|\alpha_i|>0$.
\end{remark}
From  direct computation, we find
\bel{delta_j Diff Ineq}
\begin{array}{cc}\ds
\left| \p_{\xi_1}^{\alpha_1}\prod_{i\in III} \p_{\xi_i}^{\alpha_i}   \mathfrak{S}^\sharp_j(\xi)\right|~\leq~\C_\alpha~ |\xi_1|^{-|\alpha_1|}\left[ 2^{-j}|\xi_1|\right]^{-\sum_{i\in III}|\alpha_i|}
\end{array}
\eeq
for every multi-indices $\alpha_1$ and $\alpha_i, i\in III$.

Recall $\phi_j, j\in\Z$ defined in (\ref{phi_j}).  
Consider
\bel{Delta_tj F}
\begin{array}{cc}\ds
\F^\triangle_{\t j~\left[I, \flat, \sharp\right]} f(x)~=~\int_{\R^\N} f(y)\Lambda_{\t j~\left[I, \flat, \sharp\right]} (x,y) dy,
\\\\ \ds
\Lambda_{\t j~\left[I, \flat, \sharp\right]}(x,y)~=~
\int_{\R^\N} e^{2\pi\i \big[\Phi(x,\xi)-y\cdot\xi \big]}\Delta_{\t j~\left[I, \flat, \sharp\right]}(\xi)\phi_j(\xi)\sigma(\xi)d\xi
\end{array}
\eeq
for every $j\ge0$ where $\Delta_{\t j~\left[I, \flat, \sharp\right]}$ is defined in (\ref{Delta_t j}).
\begin{remark} Suppose
$II\cup III=\emptyset$ so that $I=\{2,\ldots,n\}$. We have 
\bel{Delta, delta}
\Delta_{\t j~\left[I, \flat, \sharp\right]}(\xi)~=~ \prod_{i=2}^n \delta_{t_i}(\xi)~=~\delta_\t(\xi)
\eeq
as defined in (\ref{delta_t}).
In this case, we write 
\bel{Case Two terms}
\F^\triangle_{\t j~\left[I, \flat, \sharp\right]}~=~\F^\triangle_{\t j},\qquad \Lambda_{\t j~\left[I, \flat, \sharp\right]}~=~\Lambda_{\t j}.
\eeq
\end{remark}
Denote $\Q_r\subset\R^\N$ for $0<r<1$, so-called the {\it region of influence} associated to $\F^\triangle_{\t j}$, satisfying
\bel{Q norm}
\vol\{\Q_r\}~\leq~\C~r.
\eeq
This subset will be explicitly defined in the next section.
\v
{\bf Lemma One}~~{\it  Let $\Lambda_{\t j~\left[I, \flat, \sharp\right]}$ defined in (\ref{Delta_tj F}) for $j>0$.
Suppose $\sigma\in\S^{-{\N-1\over2}}_\rho$. We have
\bel{Result Size}
\begin{array}{cc}\ds
\int_{\R^\N} \left|\Lambda_{\t j~\left[I, \flat, \sharp\right]}(x,y) \right|dx~\leq~\C~ 2^{-j\ve \sum_{i\in II\cup III}\N_i} \prod_{i\in I} 2^{-\N_i t_i/2}
\end{array}
\eeq 
for some $\ve=\ve(\rho)>0$.

Moreover, if $I=\{2,\ldots,n\}$ and $II\cup III=\emptyset$, we have
\bel{Result Diff}
\begin{array}{ccc}\ds
\int_{\R^\N} \left|\Lambda_{\t j}(x,y)-\Lambda_{\t j}(x,0) \right|dx~\leq~\C~ 2^j|y|\prod_{i=2}^n 2^{-\N_i t_i/2}
\end{array}
\eeq 
and
\bel{Result Q}
\begin{array}{ccc}\ds
\int_{\R^\N\setminus \Q_r} \left|\Lambda_{\t j}(x,y) \right|dx~\leq~\C~
{2^{-j}\over r} \prod_{i=2}^n 2^{-\N_i t_i/2}
\end{array}
\eeq 
for  $|y|\leq r$ whenever $2^j>r^{-1}$.}

\v
Next, we prove (\ref{Result One H})  by applying (\ref{Result Size}) and (\ref{Result Diff})-(\ref{Result Q}).

Recall  {\bf Remark 1.2}. $\H^1(\R^\N)$ has an atom decomposition as characterized by   Coifman \cite{Coifman}. 
Denote $a$ to be an $\H^1$-atom: $\supp a\subset\B_r(x_o)$, $|a(x)|\leq|\B_r(x_o)|^{-1}$ and $\ds\int_{\B_r(x_o)}a(x)dx=0$ where $\B(x_o)\subset\R^\N$ is a ball centered on some $x_o$ with radius $r>0$.  

Let $\F^\triangle$ defined in (\ref{F tri}). As a  singular integral realization, we write 
\bel{F Int}
\F^\triangle f(x)~=~\int_{\R^\N} f(y) \left\{ \int_{\R^\N} e^{2\pi\i\big[\Phi(x,\xi)-y\cdot\xi\big]}\Delta(\xi)\sigma(\xi)d\xi\right\} dy.
\eeq
By changing $y\mt y-x_o$ and $x\mt x-x_o$ inside (\ref{F Int}), we can assume $\supp a\subset \B_r\doteq\B_r(0)$ without lose of the generality.

To conclude (\ref{Result One H}), it is suffice to show
\bel{H^1 est F_ta}
\begin{array}{lr}\ds
\int_{\R^\N} \left|\F^\triangle a(x)\right|dx~\leq~\C_\rho,\qquad\sigma\in\S^{-{\N-1\over 2}}_\rho.
\end{array}
\eeq 
From (\ref{F tri}) and (\ref{Delta_tj F}), we find
\bel{F tri Sum}
\F^\triangle~=~\sum_{I\cup II \cup III=\{2,\ldots,n\}} ~ \sum_{j>0}~\sum_{0\leq t_i<j/2, ~i\in I}~\F^\triangle_{\t j~\left[I,\flat,\sharp\right]}.
\eeq
\vsk
{\bf Case One}~~Suppose $II\cup III\neq\emptyset$. By using (\ref{Result Size}), we have
\bel{EST ONE}
\begin{array}{lr}\ds
\int_{\R^\N} \left|\F^\triangle_{\t j~\left[I,\flat,\sharp\right]}  a(x)\right| dx~\leq~\left\| a\right\|_{\L^1(\R^\N)} \int_{\R^\N} \left| \Lambda_{\t j~\left[I, \flat, \sharp\right]}(x,y)\right|dx 
\\\\ \ds~~~~~~~~~~~~~~~~~~~~~~~~~~~~~~~~~~
~\leq~\C~ 2^{-j\ve \sum_{i\in II\cup III}\N_i} \prod_{i\in I} 2^{-\N_i t_i/2}. 
\end{array}
\eeq
Together with {\bf Remark 3.1}, we obtain
\bel{EST ONE SUM}
\begin{array}{lr}\ds
\sum_{I\cup II \cup III=\{2,\ldots,n\}} ~ \sum_{j>0}~\sum_{0\leq t_i<j/2,~ i\in I}~\int_{\R^\N} \left|\F^\triangle_{\t j~\left[I,\flat,\sharp\right]}  a(x)\right| dx
\\\\ \ds
~\leq~\C \sum_{I\cup II \cup III=\{2,\ldots,n\}} ~ \sum_{j>0}~\sum_{0\leq t_i<j/2,~ i\in I}~2^{-j\ve \sum_{i\in II}\N_i} \prod_{i\in I} 2^{-\N_i t_i/2}\qquad\hbox{\small{ by (\ref{EST ONE})}}
\\\\ \ds
~\leq~\C_\rho.
\end{array}
\eeq

{\bf Case Two}~~Suppose $II\cup III=\emptyset$ so that $I=\{2,\ldots,n\}$. Recall  (\ref{Delta, delta})-(\ref{Case Two terms}) from {\bf Remark 3.3}.   We need to handle
\bel{F natural}
\begin{array}{cc}\ds
{^\natural}\F^\triangle ~=~  \sum_{j>0}~\sum_{0\leq t_i<j/2, ~i=2,\ldots,n}~\F^\triangle_{\t j}.
\end{array}
\eeq
Furthermore, we write
\bel{F natural Int}
\begin{array}{cc}\ds
{^\natural}\F^\triangle f(x)~=~\int_{\R^\N} e^{2\pi\i \Phi(x,\xi)} \Delta^\natural(\xi) \sigma(\xi) \Hat{f}(\xi)d\xi,
\\\\ \ds
\Delta^\natural(\xi)~=~\sum_{j>0}~\sum_{0\leq t_i<j/2, ~i=2,\ldots,n} \delta_\t(\xi)\phi_j(\xi).\end{array}
\eeq
Denote
$\ds
\left\| f\right\|_{\L^{p_1}\cdots\L^{p_n}\left(\R^{\N_1}\times\cdots\times\R^{\N_n}\right)}=\left\{\int_{\R^{\N_1}}\cdots\left\{\int_{\R^{\N_n}}\left| f(x)\right|^{p_n} dx_n\right\}^{p_{n-1}\over p_n} \cdots d x_1\right\}^{1\over p_1}
$ for
$1<p_i<\infty$,  $ i=1,2,\ldots,n$. 

Let $\sigma\in\S^{-\m}_\rho$, $0<\m\leq {\N-1\over 2}$ satisfying the differential inequality in (\ref{Class}). 
We claim
\bel{F tri L^2-L^p}
\begin{array}{cc}\ds
\left\| {^\natural}\F^\triangle f\right\|_{\L^2(\R^\N)}~\leq~\left\| f\right\|_{\L^{p_1}\cdots\L^{p_n}\left(\R^{\N_1}\times\cdots\times\R^{\N_n}\right)},
\\\\ \ds
{\rho_i\over \N_i}~=~{1\over p_i}-{1\over 2},\qquad i~=~1,2,\ldots,n.
\end{array}
\eeq
In order to prove (\ref{F tri L^2-L^p}), rewrite ${^\natural}\F^\triangle f$ in (\ref{F natural Int})
as
\bel{F natural rewrite}
\begin{array}{cc}\ds
{^\natural}\F^\triangle f(x)~=~\int_{\R^\N} e^{2\pi\i \Phi(x,\xi)} \Delta^\natural(\xi) \sigma(\xi)\prod_{i=1}^n |\xi_i|^{\rho_i} \Hat{\T f}(\xi)d\xi,
\\\\ \ds
\Hat{\T f}(\xi)~=~\Hat{f}(\xi) \prod_{i=1}^n |\xi_i|^{-\rho_i}.
\end{array}
\eeq
Observe that $|\Delta^\natural(\xi) \sigma(\xi)|\prod_{i=1}^n |\xi_i|^{\rho_i}$ is uniformly bounded for $\sigma\in\S^{-\m}_\rho$. Plancherel theorem suggests that we show $\T f$ satisfying the mixed-norm estimate in (\ref{F tri L^2-L^p}).

Note that the inverse Fourier transform of $\ds\prod_{i=1}^n |\xi_i|^{-\rho_i}$, $0<\rho_i<\N_i, ~i=1,2,\ldots,n$ equals
$\ds
 \prod_{i=1}^n {\pi^{\rho_i} \Gamma(\N_i-\rho_i)\over \pi^{\N_i-\rho_i} \Gamma(\rho_i)} \left({1\over |x_i|}\right)^{\N_i-\rho_i}$ for $ |x_i|>0,~i=1,2,\ldots,n$.
Write $x=(x_i, x'_i)\in\R^{\N_i}\times\R^{\N-\N_i}$, $i=1,2,\ldots,n$.
 By applying Hardy-Littlewood-Sobolev inequality \cite{Hardy-Littlewood}-\cite{Sobolev} respectively in every $\R^{\N_i}, i=1,2,\ldots,n$, we find
\bel{T EST mixnorm}
\begin{array}{lr}\ds
\left\| \T f\right\|_{\L^2(\R^\N)}~\leq~\C_\rho \left\{\int_{\R^\N} \left\{\int_{\R^\N} |f(u)| \prod_{i=1}^n\left({1\over |x_i-u_i|}\right)^{\N_i-\rho_i} du\right\}^2 dx\right\}^{1\over 2}
\\\\ \ds
~=~\C_\rho \left\{\iint_{\R^{\N_1}\times\R^{\N-\N_1}}\left\{\iint_{\R^{\N_1}\times\R^{\N-\N_1}} |f(u_1,u'_1)| \left({1\over|x_1-u_1|}\right)^{\N_1-\rho_1}\prod_{i=2}^n\left({1\over |x_i-u_i|}\right)^{\N_i-\rho_i} du_1du'_1\right\}^2 dx_1dx'\right\}^{1\over 2}
\\\\ \ds
~\leq~\C_\rho \left\{\int_{\R^{\N-\N_1}} \left\{\int_{\R^{\N_1}}\left\{\int_{\R^{\N-\N_1}} |f(x_1,u'_1)| \prod_{i=2}^n \left({1\over |x_i-u_i|}\right)^{\N_i-\rho_i} du'_1\right\}^{p_1} dx_1\right\}^{2\over p_1}dx'_1\right\}^{1\over 2}
\\\\ \ds
~\leq~\C_\rho \left\{\int_{\R^{\N_1}} \left\{\int_{\R^{\N-\N_1}}\left\{\int_{\R^{\N-\N_1}} |f(x_1,u'_1)| \prod_{i=2}^n \left({1\over |x_i-u_i|}\right)^{\N_i-\rho_i} du'_1\right\}^2 dx'_1\right\}^{p_1\over 2}dx_1\right\}^{1\over p_1}
\\ \ds~~~~~~~~~~~~~~~~~~~~~~~~~~~~~~~~~~~~~~~~~~~~~~~~~~~~~~~~~~~~~~~~~~~~~~~~~~~
\hbox{\small{by Minkowski integral inequality}}
\\ \ds~~~~~~~~~~~~~~~~~~~~~~~~~~~~~~~~~~~~~~~~~~~~~~~~~~~~~\vdots
\\ \ds
~\leq~\C_\rho \left\{\int_{\R^{\N_1}}\cdots\left\{\int_{\R^{\N_n}}\left| f(x)\right|^{p_n} dx_n\right\}^{p_{n-1}\over p_n} \cdots d x_1\right\}^{1\over p_1}.
\end{array}
\eeq
Now, we are ready to prove (\ref{H^1 est F_ta}) for ${^\natural}\F^\triangle$.
Suppose $\supp a\subset\B_r$ for some $r\ge1$. By using Schwartz inequality, we find
 \bel{Est r>1}
\begin{array}{lr}\ds
\int_{\B_{3r}}\left|{^\natural}\F^\triangle a(x)\right| dx
~\leq~\C~r^{\N/2}~\left\| \F a\right\|_{\L^2(\R^\N)}
\\\\ \ds~~~~~~~~~~~~~~~~~~~~~~~~~~
~\leq~\C~r^{\N/ 2}~\left\| a\right\|_{\L^2(\R^\N)}
\\\\ \ds~~~~~~~~~~~~~~~~~~~~~~~~~~
~\leq~\C~r^{\N/ 2}~r^{\N\big[-1+{1\over 2}\big]}
\\\\ \ds~~~~~~~~~~~~~~~~~~~~~~~~~~
~=~\C.
\end{array}
\eeq

On the other hand, consider $x\in\R^\N\setminus \B_{3r}$. Note that $\Phi(x,\xi)=x\cdot\xi\pm|\xi|$. We have $\nabla_{\xi_i} \left[\Phi(x,\xi)-y\cdot\xi\right]=(x-y)_i\pm {\xi_i\over |\xi|}$ for every $i=1,2,\ldots,n$. There is an $\imath\in\{1,2,\ldots,n\}$ and $\ell\in \{1,\ldots, \N_\imath\}$ such that $\left|(x-y)_{\imath\ell}\pm {\xi_{\imath\ell}\over |\xi|}\right|\sim \left|x_{\imath\ell}\right|\sim |x|$ whenever $y\in\B_r$.

Define a differential operator
\bel{D_il}
\D_{\imath\ell} f~=~(2\pi\i)^{-1} \left[ (x-y)_{\imath\ell}\pm {\xi_{\imath\ell}\over |\xi|}\right]^{-1} \partial_{\xi_{\imath\ell}} f.
\eeq
An $N$-fold integration by parts shows
\bel{Kernel imath int by parts}
\begin{array}{lr}\ds
\Lambda_{\t j}(x,y)~=~\int_{\R^\N} e^{2\pi\i \big[\Phi(x,\xi)-y\cdot\xi \big]}\left({^T}\D_{\imath\ell}\right)^N \Bigg\{ \delta_\t (\xi)\phi_j(\xi)\sigma(\xi)\Bigg\}d\xi.
\end{array}
\eeq
Let $\sigma\in\S^0_\rho$ which satisfies the differential inequality in (\ref{Class}). Moreover, recall $\delta_\t$ satisfying the differential inequality in (\ref{delta_l Diff Ineq}). 
Because $0\leq t_\imath<j/2$, we find
\bel{Diff integrand est}
\left|\left({^T}\D_{\imath\ell}\right)^N \Bigg\{ \delta_\t (\xi)\phi_j(\xi)\sigma(\xi)\Bigg\}\right|~\leq~\C_N~2^{-(j/2)N}~\left({1\over 1+|x|}\right)^N
\eeq
for every $N\ge0$.

The volume of  $\supp \delta_\t\phi_j\sigma$ is bounded by $\C 2^{j\N_1}\prod_{i=2}^n 2^{(j-t_i)\N_i}$. 
By using (\ref{Kernel imath int by parts}) and  (\ref{Diff integrand est}) with $N=2\N+1$, we have
\bel{Kernel imath Est}
\begin{array}{lr}\ds
\left| \Lambda_{\t j}(x,y)\right|~\leq~\C~ 2^{j\N_1}\prod_{i=2}^n 2^{(j-t_i)\N_i} 2^{-(j/2)(2\N+1)} \left({1\over 1+|x|}\right)^{2\N+1}
\\\\ \ds~~~~~~~~~~~~~~~~
~\leq~\C~2^{-j/2} \prod_{i=2}^n 2^{-t_i\N_i} \left({1\over 1+|x|}\right)^{2\N+1}.
\end{array}
\eeq
This size estimate of kernel further implies
 \bel{Est r>1 Outside}
\begin{array}{lr}\ds
\int_{\R^\N\setminus\B_{3r}}\left|{^\natural}\F^\triangle a(x)\right| dx~\leq~\sum_{j>0}~\sum_{0\leq t_i<j/2, ~i=2,\ldots,n}~\int_{\R^\N\setminus \B_{3r}} \left|\F^\triangle_{\t j}  a(x)\right| dx\qquad\hbox{\small{by (\ref{F natural})}}
\\\\ \ds~~~~~~~~~~~~~~~~~~~~~~~~~~~~~~~~
~\leq~\sum_{j>0}~\sum_{0\leq t_i<j/2, ~i=2,\ldots,n}~\left\| a\right\|_{\L^1(\R^\N)} \int_{\R^\N\setminus \B_{3r}} \left| \Lambda_{\t j}(x,y)\right|dx 
\\\\ \ds~~~~~~~~~~~~~~~~~~~~~~~~~~~~~~~~
~\leq~\C~\sum_{j>0}~\sum_{0\leq t_i<j/2, ~i=2,\ldots,n} 2^{-j/2} \prod_{i=2}^n 2^{-t_i\N_i}\qquad\hbox{\small{by (\ref{Kernel imath Est})}}
\\\\ \ds~~~~~~~~~~~~~~~~~~~~~~~~~~~~~~~~
~\leq~\C.
\end{array}
\eeq
Let $\supp a\subset\B_r$ for $0<r<1$. 
By applying Schwartz inequality and (\ref{Q norm}), we find
\bel{Local est}
\begin{array}{lr}\ds
\int_{\Q_r} \left|{^\natural}\F^\triangle a(x)\right|dx~\leq~\C\left\|{^\natural}\F^\triangle a\right\|_{\L^2(\R^\N)}~r^{1\over 2}. 
\end{array}
\eeq
From (\ref{Local est}), by using (\ref{F tri L^2-L^p}) with $\m={\N-1\over2}$,   we have
\bel{Local est 2}
\begin{array}{lr}\ds
\int_{\Q_r} \left|{^\natural}\F^\triangle a(x)\right|dx~\leq~\C_\rho~\left\| a\right\|_{\L^{p_1}\cdots\L^{p_n}\left(\R^{\N_1}\times\cdots\times\R^{\N_n}\right)} 
r^{1\over 2}
\\\\ \ds~~~~~~~~~~~~~~~~~~~~~~~~~~
~\leq~\C_\rho~ r^{-\N}r^{\sum_{i=1}^n{\N_i\over p_i}}~ r^{1\over 2} 
\\\\ \ds~~~~~~~~~~~~~~~~~~~~~~~~~~
~=~\C_\rho,\qquad\hbox{\small{$\sum_{i=1}^n{\N_i\over p_i}={\N\over 2}+\m={\N\over 2}+{\N-1\over 2}$}}.
\end{array}
\eeq
We are left to show
\bel{Comple est l}
\begin{array}{lr}\ds
\int_{\R^\N\setminus\Q_r} \left|{^\natural}\F^\triangle a(x)\right|dx~\leq~\C_\rho,
\qquad 
\sigma\in\S^{-{\N-1\over 2}}_\rho.
\end{array}
\eeq
Recall (\ref{F natural}). For $2^j\leq r^{-1}$, we write
\bel{F_t cancella}
\begin{array}{lr}\ds
\F^\triangle_{\t j} a(x)~=~\int_{\R^\N} a(y)\Lambda_{\t j}(x,y)dy
~=~\int_{\B_r} a(y) \left[\Lambda_{\t j}(x,y)-\Lambda_{\t j}(x,0)\right]dy
\end{array}
\eeq 
because 
$\int_{\B_r} a(y)dy=0$. 

By applying (\ref{Result Diff}), we find
\bel{Norm Est1}
\begin{array}{lr}\ds
\int_{\R^\N}\left|\F^\triangle_{\t j} a(x)\right|dx
~\leq~\int_{\B_r} |a(y)|\left\{\int_{\R^\N} \left|\Lambda_{\t j}(x,y)-\Lambda_{\t j}(x,0)\right| dx\right\} dy
\\\\ \ds~~~~~~~~~~~~~~~~~~~~~~~~~
~\leq~ \C~ 2^j r \prod_{i=2}^n 2^{-t_i\N_i/2}.
\end{array}
\eeq
By using (\ref{Norm Est1}), we have
\bel{Sum1}
\begin{array}{lr}\ds
\sum_{j\colon 2^j\leq r^{-1}}~\sum_{0\leq t_i<j/2, ~i=2,\ldots,n}~\int_{\R^\N}\left|\F^\triangle_{\t j}  a(x)\right|dx
~\leq~\C.
\end{array}
\eeq
For $2^j>r^{-1}$,  by applying (\ref{Result Q}), we find
\bel{Norm Est2}
\begin{array}{lr}\ds
\int_{\R^\N\setminus\Q_r}\left|\F^\triangle_{\t j}  a(x)\right|dx
~\leq~\int_{\B_r} |a(y)|\left\{\int_{\R^\N\setminus\Q_r} \left|\Lambda_{\t j}(x,y)\right| dx\right\} dy
\\\\ \ds~~~~~~~~~~~~~~~~~~~~~~~~~~~~~~
~\leq~\C~{2^{-j}\over r}\prod_{i=2}^n 2^{-t_i\N_i/2}.
\end{array}
\eeq
By using (\ref{Norm Est2}), we have
\bel{Sum2}
\begin{array}{lr}\ds
\sum_{j\colon 2^j> r^{-1}}~\sum_{0\leq t_i<j/2, ~i=2,\ldots,n}~\int_{\R^\N}\left|\F^\triangle_{\t j}  a(x)\right|dx
~\leq~\C.
\end{array}
\eeq
From (\ref{F natural Int}), (\ref{Sum1}) and (\ref{Sum2}), we conclude (\ref{Comple est l}).

On the other hand,  the adjoint operator of $\F^\triangle$, denoted by $(\F^\triangle)^*$ is given by (\ref{Ff}) with  $\Phi(x,\xi)=-x\cdot\xi\mp|\xi|$
and $\Delta(\xi)\sigma(\xi)$ replaced by $\bar{\Delta}(\xi)\bar{\sigma}(\xi)$.
\begin{remark}
It should be clear that $(\F^\triangle)^*\colon\H^1(\R^\N)\mt\L^1(\R^\N)$ for $\sigma\in\S^{-{\N-1\over2}}_\rho$ as (\ref{Result A H}).  
\end{remark}
Next, we prove (\ref{Result One}) by proceeding with a complex interpolation argument in analogue to 
{\bf 4.9}, chapter IX  of Stein \cite{Stein}.

Let $\sigma\in\S^{-\m}_\rho$ for $0\leq\m\leq{\N-1\over 2}$. Recall $\rho_i=\m\theta_i\ge0$ for ${\N_i-1\over \N-1}<\theta_i<{\N_i\over \N-1}$, $i=1,2,\ldots,n$. 

We consider a family of analytic operators
\bel{T_z}
\begin{array}{lr}\ds
\T_z f(x)~=~e^{(z-\vartheta)^2}\int_{\R^\N} e^{2\pi\i \Phi(x,\xi)} \Delta(\xi)\sigma(\xi)\prod_{i=1}^n \Big(1+|\xi_i|^2\Big)^{\gamma_i(z)/2} \Hat{f}(\xi)d\xi,
~~~~~
 0\leq\Re z\leq1
 \end{array}
\eeq
where 
\bel{gamma_i}
\vartheta~=~{2\m\over \N-1},\qquad \gamma_i(z)~=~\m \theta_i-z \left({\N-1\over 2}\right)\theta_i,\qquad i~=~1,2,\ldots,n.
\eeq
Observe that 
\[
e^{(z-\vartheta)^2}~=~\exp\Bigg\{ (\Re z-\vartheta)^2+2\i (\Re z-\vartheta)\Im z-(\Im z)^2\Bigg\}
\]
decays rapidly as $|\Im z|\mt\infty$. Consequently, we have
\bel{DIFF INEQ}
\left|\partial_\xi^\alpha e^{(z-\vartheta)^2} \sigma(\xi) \prod_{i=1}^n \Big(1+|\xi_i|^2\Big)^{\gamma_i(z)/2}\right|~\leq~\C_\alpha\prod_{i=1}^n \left({1\over 1+|\xi_i|}\right)^{\Re z \left({\N-1\over 2}\right)\theta_i +|\alpha_i|}
\eeq
for every multi-index $\alpha$.

By applying Plancherel theorem, we find
\bel{EST T_0}
\left\| \T_{0+\i\Im z} f\right\|_{\L^2(\R^\N)}~\leq~\C~\left\| f\right\|_{\L^2(\R^\N)}.
\eeq
Indeed,  the norm of $e^{\pm2\pi\i|\xi|}\Delta(\xi)\sigma(\xi)\prod_{i=1}^n \Big(1+|\xi_i|^2\Big)^{\gamma_i(z)/2}$   is uniformly bounded at $z=0+\i\Im z$.

On the other hand, (\ref{Result One H}) for $\F^\triangle$ implies
\bel{EST T_1}
\left\| \T_{1+\i\Im z} f\right\|_{\L^1(\R^\N)}~\leq~\C_\rho~\left\| f\right\|_{\H^1(\R^\N)}.
\eeq
Because of {\bf Remark 3.4}, 
by duality between $\H^1$ and $\BMO$ spaces, we  also have
\bel{EST* T_1}
\left\| \T_{1+\i\Im z} f\right\|_{\BMO(\R^\N)}~\leq~\C_\rho~\left\| f\right\|_{\L\infty(\R^\N)}.
\eeq
Note that 
\[
\T_\vartheta f(x)~=~ \F^\triangle f,\qquad \vartheta~=~{2\m\over \N-1}.
\]
From (\ref{EST T_0}) and (\ref{EST T_1})-(\ref{EST* T_1}), 
by applying the complex interpolation theorem set out in {\bf 5.2}, chapter IV of Stein \cite{Stein}, we obtain 
\bel{Result A tri}
  \begin{array}{cc}\ds
  \left\| \F^\triangle f\right\|_{\L^p(\R^\N)}~\leq~\C_{\rho~p}~\left\| f\right\|_{\L^p(\R^\N)},\qquad 
 \left|{1\over p}- {1\over 2}\right|~\leq~{\m\over \N-1}.
  \end{array}
 \eeq

\section{A second  dyadic decomposition}
\setcounter{equation}{0}
Let $\{2,\ldots,n\}=I\cup II\cup III$ to be defined  as (\ref{I,II}). With suitable permutations on indices, we assume $1\cup I=\{1,2,\ldots,m\}$. In particular,    $m=1$ if $I=\emptyset$ and $m=n$ for $I=\{2,\ldots,n\}$. 
Let $\M=\N_1+\N_2+\cdots+\N_m$. Moreover,  $\mathds{S}^{\M-1}$  is the unit sphere of  $\xi$-subspace  consisting $(\xi_1,\xi_2,\ldots,\xi_m)^T\in\R^{\N_1}\times\R^{\N_2}\times\cdots\times\R^{\N_m}$.

Given $j>0$,  $\left\{ \xi^\nu_j\right\}_\nu$ is a collection of points  almost equally distributed on  $\mathds{S}^{\M-1}$ having a grid length between $2^{-j/2-1}$ and $2^{-j/2}$:   {\bf (1)} $|\xi^\mu_j-\xi^\nu_j|\ge 2^{-j/2-1}$ for every $\xi^\mu_j\neq\xi^\nu_j$. {\bf (2)} For any $u\in \mathds{S}^{\M-1}$, there is a $\xi^\nu_j$ in the open set $\{\xi\in\mathds{S}^{\M-1}\colon |\xi-u|<2^{-j/2+1}\}$.

We define
\bel{phi^v_j}
\varphi^\nu_j(\xi)~=~{\ds\varphi\Bigg[2^{j/2}\left|{\xi\over |\xi|}-\xi^{\nu}_{j}\right|\Bigg]\over \ds\sum_{\nu\colon \xi^\nu_j\in\mathds{S}^{\M-1}} \varphi\Bigg[2^{j/2}\left|{\xi\over |\xi|}-\xi^{\nu}_{j}\right|\Bigg]}
\eeq
for every 
\bel{U_j}
\xi~\in~\mathcal{U}_j~=~\Bigg\{ \xi\in\R^\N\colon  \sum_{\nu\colon \xi^\nu_j\in\mathds{S}^{\M-1}} \varphi\Bigg[2^{j/2}\left|{\xi\over |\xi|}-\xi^{\nu}_{j}\right|\Bigg]\neq0\Bigg\}.
\eeq
Observe that the support of $\varphi^\nu_j$ is contained in
\bel{Gamma_j}
\Gamma^\nu_j~=~\Bigg\{\xi\in\R^\N\colon \left|{\xi\over |\xi|}-\xi^{\nu}_{j}\right|<2^{-j/2+1}\Bigg\}.
\eeq
{\bf Case One}~~Suppose $II\cap III\neq\emptyset$.
Recall   (\ref{Delta_tj F}). Given $\t$ and $j\ge0$, we have
\bel{KERNEL}
\begin{array}{cc}\ds
\Lambda_{\t j~\left[I, \flat, \sharp\right]}(x,y)~=~
\sum_{j/2\leq t_i<j, i\in II}~\Lambda_{\t j~\left[I, II, \sharp\right]}(x,y),
\\\\ \ds
\Lambda_{\t j~\left[I, II, \sharp\right]}(x,y)~=~ \int_{\R^\N} e^{2\pi\i \big[\Phi(x,\xi)-y\cdot\xi \big]}\sigma(\xi)\phi_j(\xi)\Delta_{\t j~[I,II,\sharp]}(\xi)d\xi,
\\\\ \ds
\Delta_{\t j~[I,II,\sharp]}(\xi)~=~\mathfrak{S}^\sharp_j(\xi)\prod_{i\in I} \delta_{t_i}(\xi) \prod_{i\in II} \delta_{t_i}(\xi).
\end{array}
\eeq
where $\phi_j$, $\delta_{t_i}$ and $\mathfrak{S}^\sharp_j$ are defined in (\ref{phi_j}), (\ref{delta_t}) and (\ref{Delta_t j}).
Note that 
\bel{supp flat}
\supp \sum_{j/2\leq t_i<j, i\in II}\prod_{i\in II} \delta_{t_i} ~\subset~\Bigg\{\xi\in\R^\N\colon ~~{|\xi_i|\over |\xi_1|}<2^{-j/2+1},~~i\in II\Bigg\}
\eeq
and
\bel{supp sharp}
 \supp\mathfrak{S}^\sharp_j~\subset~\Bigg\{\xi\in\R^\N\colon ~~{|\xi_i|\over |\xi_1|}<2^{-j+1},~~i\in III\Bigg\}.
\eeq
\begin{remark} As a geometric consequence of (\ref{supp flat})-(\ref{supp sharp}) , we find
\bel{supp contain}
\supp\sum_{j/2\leq t_i<j, i\in II}\Delta_{\t j~[I,II,\sharp]}~\subset~\mathcal{U}_j.
\eeq
Moreover, there are at most 
\bel{terms number}
\C~2^{j\big({\M-1\over 2}\big)}\prod_{i\in I} 2^{-t_i \N_i}
\eeq
many $\nu$'s for which $\xi^\nu_j\in\mathds{S}^{\M-1}\cap\supp\sum_{j/2\leq t_i<j, i\in II}\Delta_{\t j~[I,II,\sharp]}$.
\end{remark}

Let $\varphi^\nu_j$  defined in  (\ref{phi^v_j})-(\ref{U_j}). From (\ref{KERNEL}), we write
\bel{Lambda_tj sum}
\begin{array}{cc}\ds
\Lambda_{\t j~\left[I, II, \sharp\right]}(x,y)~=~\sum_{\nu\colon \xi^\nu_j\in\mathds{S}^{\M-1}} ~\Lambda^\nu_{\t j~\left[I, II, \sharp\right]}(x,y),
\\\\ \ds
\Lambda^\nu_{\t j~\left[I, II, \sharp\right]}(x,y)~=~
\int_{\R^\N} e^{2\pi\i \big[\Phi(x,\xi)-y\cdot\xi \big]}\sigma(\xi)\phi_j(\xi)\varphi^\nu_j(\xi)\Delta_{\t j~[I,II,\sharp]}(\xi)d\xi.
\end{array}
\eeq

{\bf Case Two} ~~Suppose $II\cup III=\emptyset$ so that $I=\{2,\ldots,n\}$. In this case, we have  $m=n$ and $\M=\N$. Moreover,  $\mathcal{U}_j=\R^\N$ in (\ref{U_j}) where the summation is taking over every $\nu$: $\xi^\nu_j\in\mathds{S}^{\N-1}$.

Recall {\bf Remark 3.2}. Let $\varphi^\nu_j$ defined in (\ref{phi^v_j}). We have
\bel{Lambda_tj I sum}
\begin{array}{cc}\ds
\Lambda_{\t j}(x,y)~=~\sum_{\nu\colon\xi^\nu_j\in\mathds{S}^{\N-1}} ~\Lambda^\nu_{\t j}(x,y),
\\\\ \ds
\Lambda^\nu_{\t j}(x,y)~=~
\int_{\R^\N} e^{2\pi\i \big[\Phi(x,\xi)-y\cdot\xi \big]}\sigma(\xi)\phi_j(\xi)\varphi^\nu_j(\xi)\delta_\t(\xi)d\xi.
\end{array}
\eeq
\begin{remark} There are at most 
\bel{terms number I}
\C~2^{j\big({\N-1\over 2}\big)}\prod_{i=2}^n 2^{-t_i \N_i}
\eeq
many $\nu$'s such that $\xi^\nu_j\in\mathds{S}^{\N-1}\cap \supp\delta_\t$.
\end{remark}
Given $\xi^\nu_j\in\mathds{S}^{\N-1}$, we consider
\bel{R xi natural}
\hbox{\bf R}^\nu_j~=~\left\{x\in\R^\N~\colon \left|\left(\nabla_\xi\Phi\right)\left(x,\xi^\nu_j\right)\cdot \xi^\nu_j\right|\leq\c 2^{-j},~\left|\left(\nabla_\xi \Phi\right)\left(x,\xi^\nu_j\right)\right|\leq\c 2^{-j/2}\right\}
\eeq
where $\c>0$  is some  constant. 
In fact,  $\Phi(x,\xi)=x\cdot\xi\pm|\xi|$ implies 
\bel{nabla Phi terms}
(\nabla_\xi \Phi)\left(x,\xi^\nu_j\right)~=~x\pm\xi^\nu_j,\qquad \left(\nabla_\xi\Phi\right)\left(x,\xi^\nu_j\right)\cdot \xi^\nu_j~=~x\cdot\xi^\nu_j\pm1.
\eeq
Recall ${^\natural}\F^\triangle$ defined in (\ref{F natural}). The region of influence associated to ${^\natural}\F^\triangle$ is defined by
  \bel{Q natural}
\Q_r~=~\bigcup_{j\colon2^{-j}\leq r} ~\bigcup_{\nu\colon\xi_j^\nu\in\mathds{S}^{\N-1}} \hbox{\bf R}^\nu_j.
\eeq
Each $\hbox{\bf R}^\nu_j$ has a volume bounded by $\C \c^\N2^{-j}2^{-j\left({\N-1\over 2}\right)}$.  There are at most a constant multiple of $2^{j\left({\N-1\over 2}\right)}$ many $\nu$'s in the collection $\left\{ \xi^\nu_j\right\}_\nu$. A direct computation shows 
\bel{Q natural norm}
\begin{array}{lr}\ds
\vol\{\Q_r\}~\leq~\C~\sum_{j\colon2^{-j}\leq r}~\sum_{\nu\colon\xi_j^\nu\in\mathds{S}^{\N-1}} 2^{-j}2^{-j\left({\N-1\over 2}\right)}
~\leq~\C~ r.
\end{array}
\eeq

\section{Proof of Lemma One}
\setcounter{equation}{0}
Let $\{2,\ldots,n\}=I\cup II\cup III$  defined in (\ref{I,II}). Without lose of the generality, we assume $1\cup I=\{1,2,\ldots,m\}$:    $m=1$ if $I=\emptyset$ and $m=n$ for $I=\{2,\ldots,n\}$. 

Recall $\M=\N_1+\N_2+\cdots+\N_m$. $\left\{ \xi^\nu_j\right\}_\nu$ is a collection of points  equally distributed on  $\mathds{S}^{\M-1}$ having a grid length between $2^{-j/2-1}$ and $2^{-j/2}$. 

Denote 
\[\eta~=~(\eta_{11},\eta^\dagger)^T~\in~\R\times\R^{\N-1},\qquad \eta^1~=~(\eta_{11},0)^T~\in\R\times\R^{\N-1},\qquad\eta^\nu_j~=~(1,0)^T~\in~\R\times\R^{\N-1}.\] 
Given $\nu$ fixed,  assert
  \bel{L block}
 \xi~=~\L_\nu\eta,\qquad \L_\nu~=~\left[\begin{array}{lr}\ds L_\nu 
\\ \ds
~~~~~ Id\end{array}\right]
\eeq
where $L_\nu$ is an $\M\times \M$-orthogonal matrix   with $\det L_\nu=1$. Moreover, 
we require 
$\xi^\nu_j=\L_\nu \eta^\nu_j$.
\v

{\bf Case One}~~Suppose $II\cup III\neq\emptyset$. Let  $\Lambda_{\t j~\left[I, II, \sharp\right]}$ defined  in (\ref{Lambda_tj sum}). Rewrite
\bel{Lambda_tj}
\begin{array}{cc}\ds
\Lambda_{\t j~\left[I, II, \sharp\right]}(x,y)~=~\sum_{\nu\colon \L_\nu\eta^\nu_j\in\mathds{S}^{\M-1}} ~\Lambda^\nu_{\t j~\left[I, II, \sharp\right]}(x,y),
\\\\ \ds
\Lambda^\nu_{\t j~\left[I, II, \sharp\right]}(x,y)~=~
\int_{\R^\N} e^{2\pi\i \big[\Phi(x,\L_\nu\eta)-y\cdot \L_\nu\eta \big]}\sigma(\L_\nu\eta)\phi_j(\L_\nu\eta)\varphi^\nu_j(\L_\nu\eta)\Delta_{\t j~[I,II,\sharp]}(\L_\nu\eta)d\eta
\end{array}
\eeq
where  $\phi_j$, $\varphi^\nu_j$  are defined in  (\ref{phi_j}) and  (\ref{phi^v_j})-(\ref{U_j}). 

Furthermore,  recall $\delta_{t_i}$, $\mathfrak{S}^\sharp_j$ defined in (\ref{delta_t}) and (\ref{Delta_t j}). We have
\bel{Delta_tj I,II}
\begin{array}{lr}\ds
\Delta_{\t j~[I,II,\sharp]}(\L_\nu \eta)~=~\mathfrak{S}^\sharp_j(\L_\nu\eta)\prod_{i\in I} \delta_{t_i}(\L_\nu \eta) \prod_{i\in II} \delta_{t_i}(\L_\nu \eta)
\\\\ \ds~~~~~~~~~~~~~~~~~~~~~~
~=~\prod_{i\in III} 1-\varphi\left[2^{-j+1}{|(\L_\nu\eta)_1|\over|\eta_i|}\right]
\prod_{i\in I} \delta_{t_i}(\L_\nu \eta) \prod_{i\in II} \delta_{t_i}(\L_\nu \eta).
\end{array}
\eeq
Consider
\bel{Phi split}
\begin{array}{cc}\ds
\Phi_\nu(x,\eta)~=~\Phi(x,\L_\nu \eta)~=~\L_\nu^T x\cdot\eta\pm|\eta|
\\\\ \ds
\Phi_\nu(x,\eta)-y\cdot \L_\nu\eta
~=~\left[\Big(\nabla_\eta\Phi_\nu\Big)\left(x,\eta_j^\nu\right)-\L_\nu^T y\right]\cdot\eta~+~\Psi(\eta)
\\\\ \ds~~~~~~~~~~~~~~~~~
~=~\L^T_\nu(x-y)\cdot\eta\pm |\eta_{11}|~+~\Psi(\eta),
\\\\ \ds
\Psi(\eta)~=~\Phi_\nu(x,\eta)-\Big(\nabla_\eta\Phi_\nu\Big)\left(x,\eta_j^\nu\right)\cdot\eta~=~\pm \left[|\eta|-|\eta_{11}|\right].
\end{array}
\eeq
From (\ref{Lambda_tj})-(\ref{Phi split}), we find
\bel{Theta}
\begin{array}{cc}\ds
\Lambda^\nu_{\t j~\left[I, II, \sharp\right]}(x,y)~=~\int_{\R^N} e^{2\pi\i\big[\L^T_\nu(x-y)\cdot\eta\pm |\eta_{11}|\big]}\Theta^\nu_{\t j~\left[I, II, \sharp\right]}(x,\eta)d\eta,
\\\\ \ds
\Theta^\nu_{\t j~\left[I, II, \sharp\right]}(x,\eta)~=~e^{2\pi\i\Psi(\eta)}\sigma(\L_\nu\eta)\phi_j(\L_\nu\eta)\varphi^\nu_j(\L_\nu\eta)\Delta_{\t j~[I,II,\sharp]}(\L_\nu\eta).
\end{array}
\eeq
Define a differential operator 
\bel{D}
\begin{array}{lr}\ds
\D~=~I+2^{2j}\partial^2_{\eta_{11}} +2^j\sum_{\ell=2}^{\N_1} \partial^2_{\eta_{1\ell}} 
+ \sum_{i\in I} 2^j\Delta_{\eta_i}
+\sum_{i\in II}2^{2(j-t_i)}\Delta_{\eta_i}+\sum_{i\in III} \Delta_{\eta_i}.
\end{array}
\eeq
Let $\sigma\in\S^{-{\N-1\over 2}}_\rho$.
Given $j>0$ and $\t=(t_1,t_2,\ldots,t_n)$, we aim to show
\bel{Theta Diff Ineq}
\left|\D^N \Theta^\nu_{\t j~\left[I, II, \sharp\right]}(x,\eta)\right|~\leq~\C_N~2^{-j\ve\sum_{i\in II\cup III}\N_i}2^{-j\left({\M-1\over 2}\right)}\prod_{i\in I}2^{t_i\big({\N_i\over 2}\big)}\prod_{i\in II}2^{-(j-t_i)\big({\N_i\over 2}\big)}
\eeq
for every $N\ge0$ and some $\ve=\ve(\rho)>0$.

{\bf 1.} Recall $\varphi^\nu_j$  defined in  (\ref{phi^v_j})-(\ref{U_j}) of which  $\supp\varphi^\nu_j\subset\Gamma^\nu_j$ is defined in (\ref{Gamma_j}). We find
 \bel{eta size}
2^{j-1}~\leq~|\eta_{11}|~\leq~2^{j+1},\qquad |\eta^\dagger|~\leq~ \C 2^{j/2}
\eeq
whenever $\L_\nu \eta\in\Gamma_j^\nu\cap \left\{2^{j-1}\leq|\eta|\leq2^{j+1}\right\}$.

 A direct computation shows
\bel{nu Diff est}
\begin{array}{cc}\ds
\left|\p_\eta^\alpha \left[\varphi_j^{\nu}\left(\L_\nu\eta\right)\right]\right|~\leq~\C_{\alpha} ~2^{|\alphaup| j/2}|\eta|^{-|\alphaup|}
\end{array}
\eeq
for every multi-index $\alpha$. Moreover, we have
\bel{deri computation 1}
\begin{array}{lr}\ds
{\p\over \p \eta_{11}}{\eta_{11}\over |\eta|}~=~{1\over |\eta|}-{\eta_{11}^2\over |\eta|^3}~=~
{1\over |\eta|^3} |\eta^\dagger|^2,
\qquad
{\p\over \p \eta_{11}}{\eta^\dagger\over |\eta|}~=~-{\eta_{11}\eta^\dagger\over |\eta|^3}.
\end{array}
\eeq
By using (\ref{deri computation 1}) and taking into account $\eta_{11}\sim2^j$ and $\eta^\dagger\lesssim2^{j/2}$, we find
\bel{deri eta norm 1}
\left| {\p\over \p \eta_{11}}{\eta\over |\eta|}\right|~\sim~2^{j/2}|\eta|^{-2}~\sim~2^{-j/2}|\eta|^{-1}.
\eeq
Note that $\partial_\eta \varphi^\nu_j(\L_\nu \eta)=0$ whenever $\left|{\eta\over |\eta|}-\eta^\nu_j\right|\leq2^{-j/2}$. 
From direct computation and by using
(\ref{deri eta norm 1}), we further have
\bel{nu Diff est 1}
\begin{array}{cc}\ds
\left|\p_{\eta_{11}}^N \left[\varphi_j^\nu \left(\L_\nu\eta\right)\right]\right|~\leq~\C_{N} ~|\eta|^{-N},\qquad N\ge0.
\end{array}
\eeq
{\bf 2.}~~Recall $\Psi(\eta)=\pm \left[|\eta|-|\eta_{11}|\right]$ as shown in (\ref{Phi split}). We claim
\bel{d Est Psi}
\left|\p_\eta^\alpha\Psi(\eta)\right|~\leq~\C_\alpha~2^{-j\alpha_{11}}2^{-(j/2)\big[|\alpha|-\alpha_{11}\big]},\qquad \L_\nu \eta\in\Gamma_j^\nu\cap \left\{2^{j-1}\leq|\eta|\leq2^{j+1}\right\}
\eeq
for every multi-index $\alpha$.

To show (\ref{d Est Psi}),  we follow the lines in
 {\bf 4.5}, chapter IX of Stein \cite{Stein}.
Note that $\Psi(\eta)$ is homogeneous of degree $1$ in $\eta$: $\Psi(r\eta)=r\Psi(\eta),r>0$.  We have
\bel{Psi est1}
\left|\partial_\eta^\alpha \Psi(\eta)\right|~\leq~\C_\alpha~|\eta|^{1-|\alpha|}
\eeq
for every multi-index $\alpha$.

Indeed, $\partial_\eta^\alpha \Psi(r\eta)=r^{|\alpha|}\partial_\gamma^\alpha\Psi(\gamma)$ for $\gamma=r\eta$. On the other hand, $\partial_\eta^\alpha \Psi(r\eta)=\partial_\eta^\alpha\left[r\Psi(\eta)\right]=r\partial_\eta^\alpha\Psi(\eta)$. Choose $r=|\eta|^{-1}$. We find $\partial_\eta^\alpha\Psi(\eta)=|\eta|^{1-|\alpha|} \partial_\gamma^\alpha\Psi(\gamma)$ where $|\gamma|=1$.

Suppose $|\alpha|-\alpha_{11}\ge2$. From (\ref{eta size}) and (\ref{Psi est1}), we have
\bel{Psi est2}
\begin{array}{lr}\ds
\left|\partial_\eta^\alpha \Psi(\eta)\right|~\leq~\C_\alpha~2^{-j\big(|\alpha|-1\big)}~=~\C_\alpha~2^{-j\alpha_{11}}2^{-j\big[|\alpha|-\alpha_{11}-1\big]} 
\\\\ \ds~~~~~~~~~~~~~~~~~~~~~~~~~~~~~~~~~~~~~~~~~~~
~\leq~\C_\alpha~2^{-j\alpha_{11}}2^{-(j/2)\big[|\alpha|-\alpha_{11}\big]}. 
\end{array}
\eeq
Next, consider $|\alpha|-\alpha_{11}=0$ or $1$. Note that $\eta^1=(\eta_{11},0)^T\in\R\times\R^{\N-1}$.
A direct computation shows
 \bel{p_eta_1 N}
 \begin{array}{lr}\ds
\left( \nabla_\eta \Psi\right)(\eta^1)~=~0,\qquad 
\Big(\partial_{\eta_{11}}^N \nabla_{\eta^\dagger}\Psi\Big)(\eta^1)
 ~=~ \Big(\nabla_{\eta^\dagger}\partial_{\eta_{11}}^N \Psi\Big)(\eta^1)~=~0,\qquad N\ge0.
 \end{array}
 \eeq 
By writing out  the Taylor expansion of $\partial_{\eta_{11}}^N \Psi(\eta)$ and $\nabla_{\eta^\dagger}\partial_{\eta_{11}}^N \Psi(\eta)$  in the $\eta^\dagger$-subspace and using $\Big(\nabla_{\eta^\dagger}\partial_{\eta_{11}}^N \Psi\Big)(\eta^1)=0$ from (\ref{p_eta_1 N}), we find
 \bel{p_eta_1 N norm}
 \partial_{\eta_{11}}^N \Psi(\eta)~=~ \O\left(|\eta^\dagger|^2|\eta|^{-N-1}\right),\qquad  
 \nabla_{\eta^\dagger} \partial_{\eta_{11}}^N \Psi(\eta)~=~ \O\left(|\eta^\dagger||\eta|^{-N-1}\right)
  \eeq
for every $N\ge0$.
  
By putting together (\ref{eta size}) and (\ref{p_eta_1 N norm}), we  have 
 \bel{p_eta_1 N Est}
  \left|\partial_{\eta_{11}}^N \Psi(\eta)\right|~\leq~\C_N  2^{-j N},\qquad \left|\nabla_{\eta^\dagger}\partial_{\eta_{11}}^N \Psi(\eta)\right|~\leq~\C_N  2^{-j N}2^{-j/2},\qquad N\ge0.
  \eeq

{\bf 3.}~Let $\phi_j$ defined in (\ref{phi_j}). We find
\bel{phi_j Diff}
\left|\partial^\alpha_\eta \phi_j(\L_\nu \eta)\right|~\leq~\C_\alpha~2^{-j|\alpha|}
\eeq
for every multi-index $\alpha$. Moreover, the support of $\phi_j(\L_\nu\eta)$ is contained in the dyadic shell: $\{2^{j-1}\leq|\eta|\leq2^{j+1}\}$.

Let $\Delta_{\t j~[I,II,\sharp]}$ defined in (\ref{Delta_tj I,II}). Suppose $\sigma\in\S^{-{\N-1\over 2}}_\rho$. We aim to show
\bel{Ineq Term}
\begin{array}{lr}\ds
\left|\p_\eta^\alpha \left[\sigma(\L_\nu\eta)\Delta_{\t j~[I,II,\sharp]}(\L_\nu\eta)\right] \right|
~\leq~\C_\alpha~2^{-j\ve\sum_{i\in II\cup III}\N_i}2^{-j\left({\M-1\over 2}\right)}\prod_{i\in I}2^{t_i\big({\N_i\over 2}\big)}\prod_{i\in II}2^{-(j-t_i)\big({\N_i\over 2}\big)}
\\\\ \ds~~~~~~~~~~~~~~~~~~~~~~~~~~~~~~~~~~~~~~~~~~~~~~~~~~~~
 2^{-j \alpha_{11}} 2^{-(j/2)\big[\sum_{\ell=2}^{\N_1}\alpha_{1\ell}+\sum_{i=2}^m |\alpha_i|\big]}\prod_{i\in II} 2^{-(j-t_i)|\alpha_i|}
\end{array}
\eeq
for every multi-index $\alpha$ whenever $\L_\nu \eta\in\Gamma_j^\nu\cap \left\{2^{j-1}\leq|\eta|\leq2^{j+1}\right\}$.

Recall (\ref{Class}). Note that $\sigma\in\S^{-{\N-1\over 2}}_\rho$ satisfies
\bel{m norm}
\begin{array}{cc}\ds
\left| \partial_\xi^\alpha\sigma(\xi)\right|~\leq~\C~ \prod_{i=1}^n \left({1\over 1+|\xi_i|}\right)^{\rho_i+|\alpha_i|},
\\\\ \ds 
{\N_i-1\over 2}<\rho_i<{\N_i\over 2},\qquad i~=~1,2,\ldots,n,\qquad \rho_1+\rho_2+\cdots+\rho_n~=~{\N-1\over 2}.
\end{array}
\eeq

From (\ref{m norm}), we further have
\bel{m norm again}
\begin{array}{cc}\ds
\left|\sigma(\xi)\right|~\leq~\C~ \left({1\over 1+|\xi_1|}\right)^{{\N_1-1\over 2}+\ve_1}\prod_{i=2}^n \left({1\over 1+|\xi_i|}\right)^{{\N_i\over 2}-\ve_i},
\\\\ \ds
0<\ve_i<{1\over 2},\qquad i=1,2,\ldots,n,\qquad \ve_1=\ve_2+\cdots+\ve_n.
\end{array}
\eeq
Suppose $\xi\in\supp \Delta_{\t j~[I,II,\sharp]}\cap\{2^{j-1}\leq|\xi|\leq2^{j+1}\}$. We find
\bel{xi_i norms}
|\xi_1|\sim2^j,\qquad |\xi_i|\sim2^{j-t_i}, ~~~~i\in I\cup II,\qquad |\xi_i|\lesssim1,~~~~i\in III.
\eeq
By bringing (\ref{xi_i norms}) to (\ref{m norm again}), we find
\bel{m norm j}
\begin{array}{lr}\ds
\left|\sigma(\xi)\right|~\leq~\C~ 2^{-j \big[{\N_1-1\over 2}+\ve_1\big]} \prod_{i\in I} 2^{-(j-t_i) \big[{\N_i\over 2}-\ve_i\big]} \prod_{i\in II} 2^{-(j-t_i) \big[{\N_i\over 2}-\ve_i\big]}
\\\\ \ds~~~~~~~~~~
~=~\C~2^{-j \big[\ve_1-\sum_{i\in I} \ve_i\big]} 2^{-j\big({\M-1\over 2}\big)} \prod_{i\in I} 2^{t_i \big[{\N_i\over 2}-\ve_i\big]} \prod_{i\in II} 2^{-(j-t_i) \big[{\N_i\over 2}-\ve_i\big]}
\\\\ \ds~~~~~~~~~~
~\leq~\C~2^{-j \big[\ve_1-\sum_{i\in I} \ve_i\big]} 2^{-j\big({\M-1\over 2}\big)} \prod_{i\in I} 2^{t_i \big({\N_i\over 2}\big)} \prod_{i\in II} 2^{-(j-t_i) \big[{\N_i\over 2}-\ve_i\big]}
\\\\ \ds~~~~~~~~~~
~=~\C~2^{-j \big[\ve_1-\sum_{i\in I} \ve_i\big]} \prod_{i\in II} 2^{(j-t_i) \ve_i} ~2^{-j\big({\M-1\over 2}\big)} \prod_{i\in I} 2^{t_i \big({\N_i\over 2}\big)} \prod_{i\in II} 2^{-(j-t_i) \big({\N_i\over 2}\big)}
\\\\ \ds~~~~~~~~~~
~\leq~\C~2^{-j \big[\ve_1-\sum_{i\in I} \ve_i\big]} \prod_{i\in II} 2^{j \ve_i/2} ~2^{-j\big({\M-1\over 2}\big)} \prod_{i\in I} 2^{t_i \big({\N_i\over 2}\big)} \prod_{i\in II} 2^{-(j-t_i) \big({\N_i\over 2}\big)}
\\\\ \ds~~~~~~~~~~
~=~\C~2^{-j \big[\ve_1-\sum_{i\in I} \ve_i-{1\over 2}\sum_{i\in II}\ve_i\big]} ~2^{-j\big({\M-1\over 2}\big)} \prod_{i\in I} 2^{t_i \big({\N_i\over 2}\big)} \prod_{i\in II} 2^{-(j-t_i) \big({\N_i\over 2}\big)}
\\\\ \ds~~~~~~~~~~
~=~2^{-j\ve\sum_{i\in II\cup III}\N_i} ~2^{-j\big({\M-1\over 2}\big)} \prod_{i\in I} 2^{t_i \big({\N_i\over 2}\big)} \prod_{i\in II} 2^{-(j-t_i) \big({\N_i\over 2}\big)}
\end{array}
\eeq
where $\ve= \big[\ve_1-\sum_{i\in I} \ve_i-{1\over 2}\sum_{i\in II}\ve_i\big]/ \sum_{i\in II\cup III}\N_i>0$.

Recall $\delta_{t_i}, i=1,2,\ldots,n$ defined in (\ref{delta_t}) satisfying the differential inequality in (\ref{delta_l Diff Ineq}). Moreover, 
$\mathfrak{S}^\sharp_j$ defined in (\ref{Delta_t j}) satisfies the differential inequality in (\ref{delta_j Diff Ineq}).

From (\ref{Delta_tj I,II}), we have
\bel{Delta_tj I,II Diff}
\begin{array}{lr}\ds
\left|\partial_\xi^\alpha \Delta_{\t j~[I,II,\sharp]}(\xi)\right|~=~\left| \partial_\xi^\alpha\left[\mathfrak{S}^\sharp_j(\xi)\prod_{i\in I} \delta_{t_i}(\xi) \prod_{i\in II} \delta_{t_i}(\xi)\right]\right|
\\\\ \ds~~~~~~~~~~~~~~~~~~~~~~~~~
~\leq~\C_\alpha~|\xi_1|^{-|\alpha_1|}\left[2^{-j}|\xi_1|\right]^{-\sum_{i\in III}|\alpha_i|} \prod_{i\in I \cup II} |\xi_i|^{-|\alpha_i|}
\end{array}
\eeq
for every multi-index $\alpha$. 

Let $\eta=(\eta_{11},\eta',\eta'')^T\in\R\times\R^{\M-1}\times\R^{\N-\M}$. Write out
\bel{matrix indices}
\begin{array}{cc}\ds
\xi_{i\ell}~=~a_{i\ell}\eta_{11}+\O(1)\cdot \eta',\qquad \ell=1,\ldots,\N_i,\qquad i=1,2,\ldots,m=1\cup I,
\\\\ \ds
\xi_{i\ell}~=~\eta_{i\ell},\qquad \ell=1,\ldots,\N_i,\qquad i\in II\cup III
\end{array}
\eeq
where $a_{i\ell}$ is the entry on the ($i+\ell$)-th row and the first column of $\L_\nu$.

Note that $0\leq t_i<j/2$ for $i\in I$. On the other hand, $|\eta'|\leq\C 2^{j/2}$ from (\ref{eta size}).
Consider $\L_\nu\eta\in\Gamma^\nu_j\cap  \supp \Delta_{\t j~[I,II,\sharp]}    \cap\{2^{j-1}\leq|\eta|\leq2^{j+1}\}$.
From (\ref{xi_i norms}) and (\ref{matrix indices}), we must have
\bel{a_il}
|a_{i\ell}|~\leq~\C 2^{-t_i},\qquad \ell=1,\ldots,\N_i,\qquad i\in I.
\eeq
By using the  estimates from (\ref{m norm}) to (\ref{a_il}) and applying the chain rule of differentiation, we obtain (\ref{Ineq Term}). By putting together the results in {\bf 1}-{\bf 3}, we conclude (\ref{Theta Diff Ineq}).
\v
Recall (\ref{Theta}). We have
\bel{Lambda-Theta}
\Lambda^\nu_{\t j~\left[I, II, \sharp\right]}(x,y)~=~\int_{\R^n} e^{2\pi\i\big[\L^T_\nu(x-y)\cdot\eta\pm |\eta_{11}|\big]}\Theta^\nu_{\t j~\left[I, II, \sharp\right]}(x,\eta)d\eta.
\eeq
Note that 
\bel{supp Theta}
\left|\supp \Theta^\nu_{\t j~\left[I, II, \sharp\right]}(x,\eta)\right|~\leq~\C~2^j 2^{j\left({\M-1\over 2}\right)} \prod_{i\in II}2^{(j-t_i)\N_i}.
\eeq
By applying (\ref{Theta Diff Ineq}), an $N$-fold integration by parts $w.r.t$ $\D$ in (\ref{D}) shows
\bel{Lambda v size}
\begin{array}{lr}\ds
\left|\Lambda^\nu_{\t j~\left[I, II, \sharp\right]}(x,y)\right|~\leq~
\\\\ \ds
\left\{1+4\pi^2 2^{2j}\left[\L_\nu^T (x-y)_{11}\pm1\right]^2
+
4\pi^2 2^j\sum_{\ell=2}^{\N_1} \left[\L_\nu^T (x-y)_{1\ell}\right]^2+4\pi^2 2^j\sum_{i=2}^{m} \left|\L_\nu^T (x-y)_{i}\right|^2\right.
\\\\ \ds
\left.~+4\pi^2 \prod_{i\in II}2^{2(j-t_i)}  \left|\L_\nu^T (x-y)_{i}\right|^2+4\pi^2 \prod_{i\in III}  \left|\L_\nu^T (x-y)_{i}\right|^2
\right\}^{-N}\int_{\R^\N} \left| \D^N \Theta^\nu_{\t j~\left[I, II, \sharp\right]}(x,\eta)\right|d\eta
\\\\ \ds
~\leq~\C_{N}~2^{-j\ve\sum_{i\in II\cup III}\N_i}2^{-j\left({\M-1\over 2}\right)}\prod_{i\in I}2^{t_i\big({\N_i\over 2}\big)}\prod_{i\in II}2^{-(j-t_i)\big({\N_i\over 2}\big)}~~2^j 2^{j\left({\M-1\over 2}\right)} \prod_{i\in II}2^{(j-t_i)\N_i}
\\\\ \ds~~~~~
\left\{1+4\pi^2 2^{2j}\left[\L_\nu^T (x-y)_{11}\pm1\right]^2
+
4\pi^2 2^j\sum_{\ell=2}^{\N_1} \left[\L_\nu^T (x-y)_{1\ell}\right]^2+4\pi^2 2^j\sum_{i=2}^{m} \left|\L_\nu^T (x-y)_{i}\right|^2\right.
\\\\ \ds~~~~~
\left.~+4\pi^2 \prod_{i\in II}2^{2(j-t_i)}  \left|\L_\nu^T (x-y)_{i}\right|^2+4\pi^2 \prod_{i\in III}  \left|\L_\nu^T (x-y)_{i}\right|^2
\right\}^{-N}.
\end{array}
\eeq
Given $\nu$ fixed and $y\in\R^\N$, denote
\bel{X}
\mathcal{X}~=~\L^T_\nu (x-y).
\eeq
From (\ref{Lambda v size}), we have
\bel{INT Size Est}
\begin{array}{lr}\ds
\int_{\R^\N} \left|\Lambda^\nu_{\t j~\left[I, II, \sharp\right]}(x,y)\right| dx~\leq~
\\\\ \ds
\C_N~2^{-j\ve\sum_{i\in II\cup III}\N_i}2^{-j\left({\M-1\over 2}\right)}\prod_{i\in I}2^{t_i\big({\N_i\over 2}\big)}\prod_{i\in II}2^{-(j-t_i)\big({\N_i\over 2}\big)}~~2^j 2^{j\left({\M-1\over 2}\right)} \prod_{i\in II}2^{(j-t_i)\N_i}
\\\\ \ds~~~~~
\int_{\R^\N} \left\{1+4\pi^2 2^{2j}\left[\mathcal{X}_{11}\pm1\right]^2
+
4\pi^2 2^j\sum_{\ell=2}^{\N_1} \left[\mathcal{X}_{1\ell}\right]^2+4\pi^2 2^j\sum_{i=2}^{m} \left|\mathcal{X}_{i}\right|^2\right.
\\\\ \ds~~~~~
\left.~+4\pi^2 \prod_{i\in II}2^{2(j-t_i)}  \left|\mathcal{X}_{i}\right|^2+4\pi^2 \prod_{i\in III}  \left|\mathcal{X}_{i}\right|^2 \right\}^{-N}d\mathcal{X}
\\\\ \ds
~\leq~\C_N~2^{-j\ve\sum_{i\in II\cup III}\N_i}2^{-j\left({\M-1\over 2}\right)}\prod_{i\in I}2^{t_i\big({\N_i\over 2}\big)}\prod_{i\in II}2^{-(j-t_i)\big({\N_i\over 2}\big)}
\\\\ \ds~~
\int_{\R^\N} \left\{1+4\pi^2 \left[\mathcal{Z}_{11}\right]^2
+
4\pi^2 \sum_{\ell=2}^{\N_1} \left[\mathcal{Z}_{1\ell}\right]^2+4\pi^2 \sum_{i=2}^{m} \left|\mathcal{Z}_{i}\right|^2+4\pi^2 \prod_{i\in II}  \left|\mathcal{Z}_{i}\right|^2+4\pi^2 \prod_{i\in III}  \left|\mathcal{Z}_{i}\right|^2 \right\}^{-N}d\mathcal{Z},
\\\\ \ds
\mathcal{Z}_{11}=2^j[\mathcal{X}_{11}\pm1],\qquad \mathcal{Z}_{1\ell}=2^{j/2}\mathcal{X}_{1\ell},~~\ell=2,\ldots,\N_1,\qquad \mathcal{Z}_i=2^{j/2}\mathcal{X}_i,~~i=2,\ldots,m
\\\\ \ds~~~~~~~~~~~~~~~~~~~~~~~~~~~~~~~~
\mathcal{Z}_i=2^{j-t_i}\mathcal{X}_i,~~i\in II,\qquad \mathcal{Z}_i=\mathcal{X}_i,~~i\in III.
\end{array}
\eeq
By choosing $N$ sufficiently large in (\ref{INT Size Est}), we find
\bel{INT Size}
\int_{\R^\N} \left|\Lambda^\nu_{\t j~\left[I, II, \sharp\right]}(x,y)\right| dx~\leq~\C~2^{-j\ve\sum_{i\in II\cup III}\N_i}2^{-j\left({\M-1\over 2}\right)}\prod_{i\in I}2^{t_i\big({\N_i\over 2}\big)}\prod_{i\in II}2^{-(j-t_i)\big({\N_i\over 2}\big)}.
\eeq
Recall (\ref{Lambda_tj}) and {\bf Remark 4.1}. We have
\bel{Int Lambda_tj Norm}
\begin{array}{lr}\ds
\int_{\R^\N}\left|\Lambda_{\t j~\left[I, II, \sharp\right]}(x,y)\right| dx~\leq~
\sum_{\nu\colon \L_\nu\eta^\nu_j\in\mathds{S}^{\M-1}} \int_{\R^\N}\left|\Lambda^\nu_{\t j~\left[I, II, \sharp\right]}(x,y)\right| dx
\\\\ \ds~~~~~~~
~\leq~\sum_{\nu\colon \L_\nu\eta^\nu_j\in\mathds{S}^{\M-1}} ~2^{-j\ve\sum_{i\in II\cup III}\N_i}2^{-j\left({\M-1\over 2}\right)}\prod_{i\in I}2^{t_i\big({\N_i\over 2}\big)}\prod_{i\in II}2^{-(j-t_i)\big({\N_i\over 2}\big)}\qquad \hbox{\small{by (\ref{INT Size})}}
\\\\ \ds~~~~~~~
~\leq~\C~2^{j\big({\M-1\over 2}\big)}\prod_{i\in I} 2^{-t_i \N_i}~2^{-j\ve\sum_{i\in II\cup III}\N_i}2^{-j\left({\M-1\over 2}\right)}\prod_{i\in I}2^{t_i\big({\N_i\over 2}\big)}\prod_{i\in II}2^{-(j-t_i)\big({\N_i\over 2}\big)}
\\ \ds~~~~~~~~~~~~~~~~~~~~~~~~~~~~~~~~~~~~~~~~~~~~~~~~~~~~~~~~~~~~~~~~~~~~~~~~~~~~~~~~~~~~~~~~~~~~~~~~~~~~~~~~~~~~~~
 \hbox{\small{by (\ref{terms number})}}
\\\\ \ds~~~~~~~
~=~\C~2^{-j\ve\sum_{i\in II\cup III}\N_i}\prod_{i\in I}2^{-t_i\big({\N_i\over 2}\big)}\prod_{i\in II}2^{-(j-t_i)\big({\N_i\over 2}\big)}.
\end{array}
\eeq
Let $\Lambda_{\t j~\left[I, \flat, \sharp\right]}$ be given as the summation in (\ref{KERNEL}). By using (\ref{Int Lambda_tj Norm}),  we find
\bel{INT KERNEL}
\begin{array}{lr}\ds
\int_{\R^\N}\left|\Lambda_{\t j~\left[I, \flat, \sharp\right]}(x,y)\right|dx~\leq~
\sum_{j/2\leq t_i<j, i\in II}\int_{\R^\N}\left|\Lambda_{\t j~\left[I, II, \sharp\right]}(x,y)\right| dx
\\\\ \ds~~~~~~~~~~~~~~~~~~~~~~~~~~~~~~~~~~~~~
~\leq~\sum_{j/2\leq t_i<j, i\in II} 2^{-j\ve\sum_{i\in II\cup III}\N_i}\prod_{i\in I}2^{-t_i\big({\N_i\over 2}\big)}\prod_{i\in II}2^{-(j-t_i)\big({\N_i\over 2}\big)}
\\\\ \ds~~~~~~~~~~~~~~~~~~~~~~~~~~~~~~~~~~~~~
~\leq~\C~2^{-j\ve\sum_{i\in II\cup III}\N_i}\prod_{i\in I}2^{-t_i\big({\N_i\over 2}\big)}.
\end{array}
\eeq
This is the estimate in (\ref{Result Size}).
\v
{\bf Case Two}~~Suppose $I=\{2,\ldots,n\}$ and $II\cup III=\emptyset$. In this case, we have $m=n, \M=\N$ and $\L_\nu=L_\nu$ as defined in (\ref{L block}).  

Recall (\ref{Lambda_tj I sum}). We have
\bel{Lambda_tj I sum rewrite}
\begin{array}{cc}\ds
\Lambda_{\t j}(x,y)~=~\sum_{\nu\colon \L_\nu\eta^\nu_j\in\mathds{S}^{\N-1}\cap\supp\delta_\t} ~\Lambda^\nu_{\t j}(x,y),
\\\\ \ds
\Lambda^\nu_{\t j}(x,y)~=~
\int_{\R^\N} e^{2\pi\i \big[\Phi(x,\L_\nu\eta)-y\cdot\L_\nu \eta \big]}\sigma(\L_\nu\eta)\phi_j(\L_\nu\eta)\varphi^\nu_j(\L_\nu\eta)\delta_\t(\L_\nu\eta)d\xi
\end{array}
\eeq
where $\delta_\t$ is defined in (\ref{delta_t}). 

Similar to (\ref{Theta}), we write
\bel{Theta I}
\begin{array}{cc}\ds
\Lambda^\nu_{\t j}(x,y)~=~\int_{\R^n} e^{2\pi\i\big[\L^T_\nu(x-y)\cdot\eta\pm |\eta_{11}|\big]}\Theta^\nu_{\t j}(x,\eta)d\eta,
\\\\ \ds
\Theta^\nu_{\t j}(x,\eta)~=~e^{2\pi\i\Psi(\eta)}\sigma(\L_\nu\eta)\phi_j(\L_\nu\eta)\varphi^\nu_j(\L_\nu\eta)\delta_\t(\L_\nu\eta)
\end{array}
\eeq
where $\Psi$ is given in (\ref{Phi split}). 

Note that $\delta_\t$ is a special case of $\Delta_{\t j~[I,II,\sharp]}$. Let $\D$ be the differential operator defined in (\ref{D}). We have
\bel{Theta Diff Ineq I}
\left|\D^N \Theta^\nu_{\t j}(x,\eta)\right|~\leq~\C_N~2^{-j\left({\N-1\over 2}\right)}\prod_{i=2}^n 2^{t_i\big({\N_i\over 2}\big)},\qquad N\ge0.
\eeq
A repeat estimate of (\ref{Lambda-Theta})-(\ref{Lambda v size}) for $II\cup III=\emptyset$ and $I=\{2,\ldots,n\}$ shows
\bel{Lambda Size I}
\begin{array}{lr}\ds
\left|\Lambda^\nu_{\t j}(x,y)\right| ~\leq~\C_N~2^{-j\left({\N-1\over 2}\right)}\prod_{i=2}^n 2^{t_i\big({\N_i\over 2}\big)}~2^j 2^{j\left({\N-1\over 2}\right)} 
\\\\ \ds
\left\{1+4\pi^2 2^{2j}\left[\L_\nu^T (x-y)_{11}\pm1\right]^2
+
4\pi^2 2^j\sum_{\ell=2}^{\N_1} \left[\L_\nu^T (x-y)_{1\ell}\right]^2+4\pi^2 2^j\sum_{i=2}^{n} \left|\L_\nu^T (x-y)_{i}\right|^2\right\}^{-N}
\end{array}
\eeq
for every $N\ge0$.

In analogue to (\ref{X})-(\ref{INT Size}), we obtain
\bel{INT Size I}
\int_{\R^\N} \left|\Lambda^\nu_{\t j}(x,y)\right| dx~\leq~\C_N~2^{-j\left({\N-1\over 2}\right)}\prod_{i=2}^n 2^{t_i\big({\N_i\over 2}\big)}.
\eeq
Recall {\bf Remark 4.2}. From (\ref{Lambda_tj I sum rewrite}) and (\ref{INT Size I}), we have
\bel{INT Lambda_tj I sum}
\begin{array}{lr}\ds
\int_{\R^\N}\left|\Lambda_{\t j}(x,y)\right|dx~\leq~\C\sum_{\nu\colon \L_\nu\eta^\nu_j\in\mathds{S}^{\N-1}\cap\supp\delta_\t} 
2^{-j\left({\N-1\over 2}\right)}\prod_{i=2}^n 2^{t_i\big({\N_i\over 2}\big)}
\\\\ \ds~~~~~~~~~~~~~~~~~~~~~~~~~~~~
~\leq~\C ~2^{j\left({\N-1\over 2}\right)}\prod_{i=2}^n 2^{-t_i\N_i}~2^{-j\left({\N-1\over 2}\right)}\prod_{i=2}^n 2^{t_i\big({\N_i\over 2}\big)}~=~\C\prod_{i=2}^n2^{-t_i\big({\N_i\over 2}\big)} .
\end{array}
\eeq
By carrying out the same argument as (\ref{Lambda_tj I sum rewrite})-(\ref{INT Lambda_tj I sum}), we find
\bel{INT Lambda_tj I sum y}
\begin{array}{lr}\ds
\int_{\R^\N}\left|\nabla_y\Lambda_{\t j}(x,y)\right|dx~\leq~\C~2^j\prod_{i=2}^n2^{-t_i\big({\N_i\over 2}\big)}.
\end{array}
\eeq
This further implies 
\bel{Est2 >} 
\int_{\R^\N} \left|\Lambda_{\t j}(x,y)-\Lambda_{\t j}(x,0)\right| dx~\leq~\C_{\sigma~\Phi}~2^j|y|~\prod_{i=2}^n2^{-t_i\big({\N_i\over 2}\big)}.
  \eeq

Recall the region of influence $\Q_r$, $0<r<1$ defined in (\ref{R xi natural})-(\ref{Q natural}). 
Let $k\in\Z$ and $2^{k-1}< r^{-1}<2^k$. 
  Given $\xi^\nu_j$, $j\ge k$, there exists a $\xi^\mu_k\in \left\{\xi_k^\mu\right\}_\mu\subset\mathds{S}^{\N-1}$ such that $|\xi_j^\nu-\xi_k^\mu|\leq 2^{-k/2}$. 
 Suppose $x\in\R^\N\setminus\Q_r$, we must have
 \bel{Condition k}
 \left|\left(\nabla_\xi\Phi\right)\left(x,\xi^\mu_k\right)\cdot\xi^\mu_k\right|~=~\left|x\cdot\xi^\mu_k\pm1\right|~>~\c 2^{-k}\qquad\hbox{or}\qquad 
  \left|\left(\nabla_\xi \Phi\right)\left(x,\xi^\mu_k\right)\right|~=~\left|x\pm\xi^\mu_k\right|~>~\c 2^{-k/2}.
  \eeq 
 Denote $\xi^\mu_k=\L_\nu \eta^\mu_k$. Because $|\eta^\nu_j|=|(\eta^\nu_j)_{11}|=1$, we find 
 \bel{eta^mu_k}
 \begin{array}{cc}\ds
 1-2^{-k/2}~<~\left|\left(\eta^\mu_k\right)_{11}\right|~\leq~1\qquad\hbox{and}\qquad \left|\left(\eta^\mu_k\right)_{1\ell}\right|~<~2^{-k/2},\qquad \ell=2,\ldots,\N_1,
 \\\\ \ds
 \left|\left(\eta^\mu_k\right)_{i\ell}\right|~<~2^{-k/2},\qquad \ell=1,2,\ldots,\N_i,\qquad i=2,\ldots,n.
 \end{array}
 \eeq 
 Moreover, 
 $\left(\nabla_\xi \Phi\right)\left(x,\xi^\mu_k\right)\cdot\xi^\mu_k=x\cdot\xi^\mu_k\pm1=x\cdot\L_\nu\eta^\mu_k\pm1 =\left(\nabla_\eta \Phi_\nu\right) \left(x, \eta^\mu_k\right)\cdot\eta^\mu_k$.
 The first inequality in (\ref{Condition k}) together with (\ref{eta^mu_k}) imply
 \bel{condition k}
 \begin{array}{lr}\ds
\Big|\left(\partial_{\eta_{11}}\Phi_\nu\right)\left(x,\eta_k^\mu\right)\Big|~>~\left[1-{1\over \sqrt{2}}\right]\c 2^{-k}\qquad\hbox{or}
\\\\ \ds
\sum_{\ell=2}^{\N_1} \Big|\left(\partial_{\eta_{1\ell}}\Phi_\nu\right)\left(x,\eta_k^\mu\right)\Big|+
\sum_{i=2}^n\sum_{\ell=1}^{\N_i} \Big|\left(\partial_{\eta_{i\ell}}\Phi_\nu\right)\left(x,\eta_k^\mu\right)\Big|~>~ \c 2^{-k/2}.
\end{array}
\eeq
Consider
\bel{nabla Phi diff}
\begin{array}{lr}\ds
\left(\partial_{\eta_{11}}\Phi_\nu\right)\left(x,\eta_k^\mu\right)-\left(\partial_{\eta_{11}}\Phi_\nu\right)\left(x,\eta_j^\nu\right)
~=~ \left(\nabla_\eta\partial_{\eta_{11}}\Phi_\nu\right)\left(x,\eta^\nu_j\right)\cdot\left( \eta_k^\mu- \eta_j^\nu\right)+\O(1)  \left|\eta_k^\mu- \eta_j^\nu\right|^2
 \\\\ \ds~~~~~~~~~~~~~~~~~~~~~~~~~~~~~~~~~~~~~~~~~~~~~~~~~~~~~~~
 ~=~ \left(\partial_{\eta_{11}}\nabla_\eta\Phi_\nu\right)\left(x, \eta^\nu_j\right)\cdot\left( \eta_k^\mu- \eta_j^\nu\right)+\O(1)  \left|\eta_k^\mu- \eta_j^\nu\right|^2. 
\end{array}
\eeq
Note that $\nabla_\eta\Phi_\nu(x,\eta)=\L_\nu^T x\pm {\eta\over |\eta|}$.  By using (\ref{deri computation 1}),
we find $\left(\partial_{\eta_{11}}\nabla_\eta\Phi_\nu\right)\left(x, \eta^\nu_j\right)=0$ which implies
\bel{condition1}
\begin{array}{lr}\ds
\left|\left(\partial_{\eta_{11}}\Phi_\nu\right)\left(x,\eta_k^\mu\right)-\left(\partial_{\eta_{11}}\Phi_\nu\right)\left(x,\eta_j^\nu\right)\right|~\leq~\C ~\left(2^{-k/2} \right)^2
~=~\C ~2^{-k}.
\end{array}
\eeq
On the other hand, the mean value theorem implies
\bel{condition2}
\left| \left(\nabla_\eta\Phi_\nu\right)\left(x,\eta_k^\mu\right)-\left(\nabla_\eta\Phi_\nu\right)\left(x,\eta_j^\nu\right)\right|~\leq~\C~ 2^{-k/2}.
\eeq
Suppose $|y|\leq r, 0<r<1$. 
From (\ref{condition k}), (\ref{condition1}) and (\ref{condition2}), by using the triangle inequality and taking into account for $\left|\L_\nu^T y\right|\leq 2^{-k}$, we obtain
\bel{condition jk}
\begin{array}{lr}\ds
4\pi^2 2^{2j}\left[\left(\nabla_\eta\Phi_\nu\right)\left(x,\eta_j^\nu\right)-\L_\nu^T y\right]_{11}^2+
4\pi^2 2^{j}\sum_{\ell=2}^{\N_1} \left[\left(\nabla_\eta\Phi_\nu\right)\left(x,\eta_j^\nu\right)-\L_\nu^T y\right]_{1\ell}^2
\\\\ \ds
+4\pi^2 2^{j}\sum_{i=2}^n \left|\left(\nabla_\eta\Phi_\nu\right)\left(x,\eta_j^\nu\right)-\L_\nu^T y\right|_i^2
\\\\ \ds
~=~4\pi^2 2^{2j}\left[\L_\nu^T (x-y)_{11}\pm1\right]^2
+
4\pi^2 2^j\sum_{\ell=2}^{\N_1} \left[\L_\nu^T (x-y)_{1\ell}\right]^2+4\pi^2 2^j\sum_{i=2}^{n} \left|\L_\nu^T (x-y)_{i}\right|^2
~\ge~2^{j-k}
\end{array}
\eeq
provided that $\c>0$ inside (\ref{condition k}) is large enough.

Furthermore, $(\nabla_\xi\Phi)(x,\xi^\mu_k)=x\pm \xi^\mu_k =\L_\nu \L_\nu^T x\pm \L_\nu \eta^\mu_k=\L_\nu(\nabla_\eta\Phi_\nu)(x,\eta^\mu_k)$.
The second inequality in (\ref{Condition k}) implies
$\left| \left(\nabla_\eta\Phi_\nu\right) \left(x,\eta_k^\mu\right)\right|>\c2^{-k/2}$.
Together with (\ref{condition2}), by using the triangle inequality and taking into account for $\left|\L_\nu^T y\right|\leq 2^{-k}$, we find (\ref{condition jk}) again if $\c>0$ is sufficiently large.

Consider $j\ge k$ and $x\in\R^\N\setminus\Q_r$. 
From (\ref{Lambda Size I}) and (\ref{condition jk}), we have 
\bel{Lambda Size I 1-N}
\begin{array}{lr}\ds
\left|\Lambda^\nu_{\t j}(x,y)\right| ~\leq~\C_N~2^{-j\left({\N-1\over 2}\right)}\prod_{i=2}^n 2^{t_i\big({\N_i\over 2}\big)} ~2^j 2^{j\left({\N-1\over 2}\right)} 
\\\\ \ds
\left\{1+4\pi^2 2^{2j}\left[\L_\nu^T (x-y)_{11}\pm1\right]^2
+
4\pi^2 2^j\sum_{\ell=2}^{\N_1} \left[\L_\nu^T (x-y)_{1\ell}\right]^2+4\pi^2 2^j\sum_{i=2}^{n} \left|\L_\nu^T (x-y)_{i}\right|^2\right\}^{-N}
\\\\ \ds
~\leq~\C_N~2^{-j\left({\N-1\over 2}\right)}\prod_{i=2}^n 2^{t_i\big({\N_i\over 2}\big)}~2^{-(j-k)}~2^j 2^{j\left({\N-1\over 2}\right)} 
\\\\ \ds
\left\{1+4\pi^2 2^{2j}\left[\L_\nu^T (x-y)_{11}\pm1\right]^2
+
4\pi^2 2^j\sum_{\ell=2}^{\N_1} \left[\L_\nu^T (x-y)_{1\ell}\right]^2+4\pi^2 2^j\sum_{i=2}^{n} \left|\L_\nu^T (x-y)_{i}\right|^2\right\}^{1-N}.
\end{array}
\eeq
By using (\ref{Lambda Size I 1-N}) instead of (\ref{Lambda Size I}), repeat the estimates in (\ref{INT Size I})-(\ref{INT Lambda_tj I sum}).
 Finally, we arrive at
\bel{Result Q l}
\int_{\R^\N\setminus \Q_r} \left|\Lambda^\nu_{\t j}(x,y)\right|dx~\leq~\C~
{2^{-j}\over r} ~ \prod_{i=2}^n 2^{-t_i\big({\N_i\over 2}\big)}.
\eeq

\appendix
 \section{Some estimates related to Bessel functions}
 \setcounter{equation}{0}
First, we list a number of well known results for Bessel functions. More discussion  can be found in  the book of Watson \cite{Watson}.
\v

$\bullet$ For  $a>-{1\over 2}, b\in\R$ and $\rho>0$,  a Bessel function $\J_{a+\i b}$ has an integral formula 
\bel{Bessel}
\J_{a+\i b}(\rho)~=~{(\rho/2)^{a+\i b}\over \pi^{1\over 2}\Gamma\left(a+{1\over 2}+\i b\right)} \int_{-1}^1 e^{\i\rho s } (1-s^2)^{a-{1\over 2}+\i b} ds.
\eeq
$\bullet$ For every $a>-{1\over 2}, b\in\R$ and $\rho>0$, we have 
\bel{J asymptotic}
\begin{array}{lr}\ds
\J_{a+\i b}(\rho) ~\sim~\left({\pi\rho\over 2}\right)^{-{1\over 2}} \cos\left[\rho-(a+\i b) {\pi\over 2}-{\pi\over 4}\right]
\\\\ \ds~~~~~~~~~~~~
~+~\left({\pi\rho\over 2}\right)^{-{1\over 2}} \sum_{k=1}^\infty \cos\left[\rho-(a+\i b) {\pi\over 2}-{\pi\over 4}\right] \a_k \rho^{-2k}+\sin\left[\rho-(a+\i b) {\pi\over 2}-{\pi\over 4}\right]~ \b_k \rho^{-2k+1},
\\\\ \ds
\a_k~=~(-1)^k[a+\i b, 2k]2^{-2k},\qquad \b_k~=~(-1)^{k+1}[a+\i b, 2k-1]2^{-2k+1},
\\\\ \ds
[a+\i b,m]~=~{\Gamma\left({1\over 2}+a+\i b+m\right)\over m!\Gamma\left({1\over 2}+a+\i b-m\right)},\qquad  m~=~0,1,2,\ldots
\end{array}
\eeq
in the sense of that
\bel{J O}
\begin{array}{lr}\ds
\left({d\over d\rho}\right)^\ell\left[ \J_{a+\i b}(\rho) ~-~\left({\pi\rho\over 2}\right)^{-{1\over 2}} \cos\left[\rho-(a+\i b) {\pi\over 2}-{\pi\over 4}\right]\right.
\\\\ \ds~~~
\left.~-~\left({\pi\rho\over 2}\right)^{-{1\over 2}}\sum_{k=1}^N  \cos\left[\rho-(a+\i b) {\pi\over 2}-{\pi\over 4}\right] \a_k \rho^{-2k}+\sin\left[\rho-(a+\i b) {\pi\over 2}-{\pi\over 4}\right]~ \b_k \rho^{-2k+1}\right]
\\\\ \ds
~=~\O\left(\rho^{-2N-{1\over 2}}\right)\qquad N\ge0,\qquad \ell\ge0
\end{array}
\eeq
 as $\rho\mt\infty$.

$\bullet$  For every $a,b\in\R$ and $\rho>0$, we have
\bel{J identity}
\J_{a-1+\i b} (\rho)~=~2\left[{a+\i b\over \rho}\right] \J_{a+\i b}(\rho)-\J_{a+1+\i b}(\rho).
\eeq
By using (\ref{Bessel}) and (\ref{J asymptotic})-(\ref{J O}) together with (\ref{J identity}), we find the norm estimate as follows.

$\bullet$ For  every $a, b\in\R$ and $\rho>0$, 
\bel{J norm}
\left| {1\over\rho^{a+\i b}}~\J_{a+\i b}(\rho)\right|~\leq~\C_a~\left({1\over 1+\rho}\right)^{{1\over 2}+a}~e^{\hbox{\small{\bf c}}|b|}.
\eeq

Let $z\in\Cx$.  
$\Omega^z$  is a distribution defined by analytic continuation from
\bel{Omega^z}
\begin{array}{ccc}\ds
\Re z<1,\qquad\Omega^z(x)~=~\pi^{-z}\Gamma^{-1}\left(1-z\right)  \left({1\over 1-|x|^2}\right)^{z}_+.
\end{array}
\eeq
$\bullet$  $\Omega^z$ can be equivalently defined by 
 \bel{Omega^z Transform} 
\begin{array}{lr}\ds
\Hat{\Omega}^z(\xi)~=~\left({1\over|\xi|}\right)^{{\N\over 2}-z} \J_{{\N\over 2}-z}\Big(2\pi|\xi|\Big),\qquad z\in\Cx.
\end{array}
\eeq
\begin{remark} $\Hat{\Omega}^z$ in (\ref{Omega^z Transform}) is analytic for $z\in\Cx$. 
\end{remark}
By the principal of analytic continuation, the two definitions  in (\ref{Omega^z}) and (\ref{Omega^z Transform}) are equivalent.
Because (\ref{Omega^z Transform}) plays an fundamental role in our analysis, we derive this formula in below. 

Denote
\bel{area}
\omega_{\N-2}~=~2 \pi^{\N-1\over 2} \Gamma^{-1}\left({\N-1\over 2}\right)
\eeq
which is the area of $\mathds{S}^{\N-2}$.

From direct computation, we have
\bel{Omega^lambda Fourier transform}
\begin{array}{lr}\ds
\Hat{\Omega}^z(\xi)~=~\pi^{-z}  \Gamma^{-1}\left(1-z\right)\int_{|x|<1} e^{-2\pi\i x\cdot\xi} \left({1\over 1-|x|^2}\right)^z dx
\\\\ \ds~~~~~~~~~
~=~\pi^{-z}  \Gamma^{-1}\left(1-z\right)\int_0^\pi\left\{\int_0^1 e^{-2\pi\i |\xi|r\cos\vartheta} (1-r^2)^{-z} r^{\N-1}dr\right\} \omega_{\N-2} \sin^{\N-3}\vartheta d\vartheta
\\\\ \ds~~~~~~~~~
~=~\omega_{\N-2}\pi^{-z}  \Gamma^{-1}\left(1-z\right)\int_{-1}^1\left\{\int_0^1 e^{2\pi\i|\xi|rs} (1-r^2)^{-z} r^{\N-1}dr\right\} (1-s^2)^{\N-3\over 2} ds~~~~ \hbox{\small{($-s=\cos\vartheta$)}}
\\\\ \ds~~~~~~~~~
~=~\omega_{\N-2}\pi^{-z}  \Gamma^{-1}\left(1-z\right)\int_0^1\left\{\int_{-1}^1 e^{2\pi\i|\xi|rs} (1-s^2)^{\N-3\over 2} ds\right\}(1-r^2)^{-z} r^{\N-1}dr
\\\\ \ds~~~~~~~~~
~=~\omega_{\N-2}\pi^{-z}  \Gamma^{-1}\left(1-z\right)\int_0^1\left\{\int_{-1}^1 \cos(2\pi|\xi|rs) (1-s^2)^{\N-3\over 2} ds\right\}(1-r^2)^{-z} r^{\N-1}dr.
\end{array}
\eeq
Recall the Beta function identity:
\bel{Beta}
{\Gamma(z)\Gamma(w)\over \Gamma(z+w)}~=~\int_0^1 r^{z-1}(1-r)^{w-1}dr
\eeq
for every $\Re z>0$ and $\Re\ w>0$.

By writing out the Taylor expansion of the cosine function inside (\ref{Omega^lambda Fourier transform}), we find
\bel{Omega^lambda Fourier transform Sum}
\begin{array}{lr}\ds
\omega_{\N-2}\pi^{-z}  \Gamma^{-1}\left(1-z\right)\int_0^1\left\{\int_{-1}^1 \cos(2\pi|\xi|rs) (1-s^2)^{\N-3\over 2} ds\right\}(1-r^2)^{-z} r^{\N-1}dr
\\\\ \ds
=~\omega_{\N-2}\pi^{-z}  \Gamma^{-1}\left(1-z\right)
\\ \ds~~~~~
\sum_{k=0}^\infty (-1)^k {(2\pi|\xi|)^{2k}\over (2k)!}\left\{\int_{-1}^1 s^{2k} (1-s^2)^{\N-3\over 2} ds\right\}\left\{\int_0^1 r^{2k+\N-1}(1-r^2)^{-z} dr\right\}
\\\\ \ds
=~{1\over 2}\omega_{\N-2}\pi^{-z}  \Gamma^{-1}\left(1-z\right)
\\ \ds~~~~~
\sum_{k=0}^\infty (-1)^k {(2\pi|\xi|)^{2k}\over (2k)!}\left\{\int_0^1 t^{k+{1\over 2}-1} (1-t)^{{\N-1\over 2}-1} dt\right\}\left\{\int_0^1\rho^{k+{\N\over 2}-1}(1-\rho)^{1-z-1} d\rho\right\}
\\\\ \ds
=~{1\over 2}\omega_{\N-2}\pi^{-z}  \Gamma^{-1}\left(1-z\right)\sum_{k=0}^\infty (-1)^k {(2\pi|\xi|)^{2k}\over (2k)!}~{\Gamma\left(k+{1\over 2}\right)\Gamma\left({\N-1\over 2}\right)\over\Gamma\left(k+{\N\over 2}\right)}~{ \Gamma\left(k+{\N\over 2}\right) \Gamma\left(1-z\right)         \over    \Gamma\left(k+{\N\over 2}+1-z\right)} \qquad
 \hbox{\small{by (\ref{Beta})}}
\\\\ \ds
=~ \pi^{{\N-1\over 2}-z}\sum_{k=0}^\infty (-1)^k {(2\pi|\xi|)^{2k}\over (2k)!}~{\Gamma\left(k+{1\over 2}\right)        \over    \Gamma\left(k+{\N\over 2}+1-z\right)},\qquad\hbox{\small{by (\ref{area})}}.
\end{array}
\eeq
On the other hand, we have
\bel{Omega cos}
\begin{array}{lr}\ds
\left({1\over|\xi|}\right)^{{\N\over 2}-z} \J_{{\N\over 2}-z}\Big(2\pi|\xi|\Big)
~=~\hbox{\small{$\pi^{{\N-1\over 2}-z}  \Gamma^{-1}\left({\N+1\over 2}-z\right)$}}\int_{-1}^1 e^{2\pi\i  |\xi| s} (1-s^2)^{{\N-1\over 2}-z} ds\qquad\hbox{\small{by (\ref{Bessel})}}
\\\\ \ds~~~~~~~~~~~~~~~~~~~~~~~~~~~~~~~~~~
~=~\hbox{\small{$\pi^{{\N-1\over 2}-z}  \Gamma^{-1}\left({\N+1\over 2}-z\right)$}}\int_{-1}^1 \cos\left(2\pi  |\xi| s\right) (1-s^2)^{{\N-1\over 2}-z} ds
\\\\ \ds~~~~~~~~~~~~~~~~~~~~~~~~~~~~~~~~~~
~=~\hbox{\small{$\pi^{{\N-1\over 2}-z}  \Gamma^{-1}\left({\N+1\over 2}-z\right)$}} \sum_{k=0}^\infty (-1)^k {(2\pi|\xi|)^{2k}\over (2k)!}\int_{-1}^1 s^{2k} (1-s^2)^{{\N-1\over 2}-z} ds
\\\\ \ds~~~~~~~~~~~~~~~~~~~~~~~~~~~~~~~~~~
~=~\hbox{\small{$\pi^{{\N-1\over 2}-z}  \Gamma^{-1}\left({\N+1\over 2}-z\right)$}} \sum_{k=0}^\infty (-1)^k {(2\pi|\xi|)^{2k}\over (2k)!}\int_0^1 \rho^{k+{1\over 2}-1} (1-\rho)^{{\N+1\over 2}-z-1} d\rho
\\\\ \ds~~~~~~~~~~~~~~~~~~~~~~~~~~~~~~~~~~
~=~\pi^{{\N-1\over 2}-z} \sum_{k=0}^\infty (-1)^k {(2\pi|\xi|)^{2k}\over (2k)!}~{\Gamma\left(k+{1\over 2}\right)    \over    \Gamma\left(k+{\N\over 2}+1-z\right)}\qquad\hbox{\small{by (\ref{Beta}).}}
\end{array}
\eeq

{\small Westlake University, Hangzhou, 310014, China}\\
{\small email: wangzipeng@westlake.edu.cn}
\v
 
 {\small School of Mathematical Sciences, Zhejiang University, Hangzhou, 310058, China}\\

 {\small Fudan University, Shanghai, 200433, China}\\


\begin{thebibliography}{100}

\bibitem{Strichartz}{\small R. Strichartz, {\it Convolutions with kernels having singularity on a sphere}, Transaction of the American Mathematical Society {\bf 148}: 461-471, 1970}.







{\small \bibitem{Seeger-Sogge-Stein}A.~Seeger,~C.~D.~Sogge and E.~M.~Stein, {\it Regularity Properties of Fourier Integral Operators}, Annals of Mathematics, {\bf 134}: 231-251, 1991.}






{\small \bibitem{Miyachi}  A. ~Miyachi, {\it On some singular Fourier multipliers}, J. Fac. Sci. Tokyo, Sci. IA {\bf 28}: 267-315, 1981.}



\bibitem{Fefferman'}{\small C.~Fefferman, {\it A note on spherical summation multipliers},  Israel J. Math {\bf 15}: 44-52, 1973.} 







{\small \bibitem{Fefferman-Stein}C.~Fefferman and E.~M.~Stein, {\it $\H^p$ Spaces of Several Variables}, Acta Mathematica, {\bf 129}: 137-193, 1972. }

    {\small\bibitem{Coifman}R.~R.~Coifman, {\it A real variable characterization of $\H^p$}, Studia Mathematica {\bf 51}: 269-274, 1974.}



{\small \bibitem{Stein} E.~M.~Stein, {\it
Harmonic Analysis: Real-Variable Methods, Orthogonality and Oscillatory Integrals},
 Princeton University Press, 1993.}
 
 
 
 
 
 
 
 
 
 
{\small \bibitem{R.Fefferman-Stein}R.~Fefferman and E.~M.~Stein, {\it Singular Integrals on Product Spaces},
Advances in Mathematics {\bf 45}:117-143, 1982.}

 \bibitem{Rubin}{\small B.~S.~Rubin, {\it Multiplier Operators connected with the Cauchy Problem for the Wave Equation. Difference Regularization}, Math. USSR Sbornik {\bf 68}: No.2, 391-416, 1991.}
 
 
 {\small \bibitem{Journe'} J.~L.~Journ\'{e}, {\it Calder\'{o}n-Zygmund Operators on Product Spaces}, Revista Matematica Iberoamericana {\bf 1}: no.3, 55-91, 1985. }




{\small \bibitem{R.Fefferman}R.~Fefferman, {\it Harmonic Analysis on Product Spaces}, Annals of Mathematics {\bf 126}: 109-130, 1987.}
 
 
 
 
 
 


\bibitem{Cordoba-Fefferman}{\small A.~Cordoba and R.~Fefferman, {\it A geometric Proof of the Strong Maximal Theorem}, Annals of Mathematics {\bf 102}: no.1, 95-100, 1975.}



{\small \bibitem{Muller-Ricci-Stein}D. M\"{u}ller,~~F.~Ricci,~~E.~M.~Stein,  
{\it Marcinkiewicz Multipliers and Multi-parameter structures on Heisenberg (-type) group, I}, 
Inventiones Mathematicae {\bf 119}: 199-233, 1995.}




{\small\bibitem{Carleson}L.~Carleson, {\it A counterexample for measures bounded on $\H^p$ for the bi-disc}, Mittag Keffker Report, no.7, 1974.}




{\small \bibitem{Malliavin} M.~P.~Malliavin and P.~Malliavin, {\it Int\'{e}grales de Lusin-Calder\'{o}n pour les fonctions biharmoniques}, Bull. Sci. Math. {\bf 101}: 357-384, 1977.}

{\small\bibitem{Gundy-Stein}R.~F.~Gundy and E.~M. Stein, {\it $H^p$ theory for the poly disc}, Proc. Nat. Acad. Sci. USA {\bf 76}: 1026-1029, 1979.}






{\small \bibitem{Hormander}L. H\"{o}rmander, {\it Fourier integral operators I}, Acta Mathematica {\bf 127}: 79-183, 1971.}




{\small \bibitem{Duistermaat-Hormander}J. J. Duistermaat and L. H\"{o}rmander, {\it Fourier integral operators II}, Acta Mathematica {\bf 122}: 183-269, 1972.}



{\small \bibitem{Colin-Frisch}  Y. Colin de Verdi\'{e}re and M. Frisch, {\it R\'{e}gularit\'{e} Lipschitzienne et solutions de l'\'{e}quation des ondes sur une vari\'{e}t\'{e} Riemannienne compacte}, Ann. Scient. Ecole Norm. Sup. {\bf 9}: 539-565, 1976.}











{\small \bibitem{Brenner}P. Brenner, {\it $\L^p\mt\L^{p'}$-estimates for Fourier integral operators related to hyperbolic equations}, Math. Z. {\bf 152}: 273-286, 1977.}

{\small \bibitem{Peral}J. Peral, {\it $\L^p$-estimates for the wave equation}, Journal of Functional Analysis {\bf 36}: 114-145, 1980.}

{\small \bibitem{Beals}M. Beals, {\it $\L^p$-Boundedness of Fourier Integral Operators}, Mem. Amer. Math. Soc. {\bf 264}: 1982.}




















































{\small \bibitem{Hardy-Littlewood}
G.~H.~Hardy and J.~E.~Littlewood, {\it Some Properties of Fractional Integrals}, 
Mathematische Zeitschrift {\bf 27}: 565-606, 1928.}






{\small \bibitem{Sobolev}
S.~L.~Sobolev, {\it On a Theorem of Functional Analysis}, Matematicheskii Sbornik {\bf 46}: 471-497, 1938.}











\bibitem{Watson}{\small G.~N.~Watson, {\it Theory of Bessel Functions}, Cambridge University Press, Cambridge, 1944}.

{\small\bibitem{Seeger}
A. Seeger, {Degenerate Fourier integral operators in the plane}, Duke Math. J. {\bf 71}: no.3, 685-745, 1993.
}
	


{\small\bibitem{Lopez-Rule-Staubach}
S.	Rodriguez-L\'{o}pez, D.  Rule and W. Staubach, {\it  A Seeger-Sogge-Stein theorem for bilinear Fourier integral operators}, Advances in Mathematics {\bf264}: 1-54, 2014.}







{\small\bibitem{Kato-Miyachi-Tomita} 
 T.  Kato,  A. Miyachi and N. Tomita,  {\it Estimates for a certain bilinear Fourier integral operator}, Journal of Pseudo-Differential Operators and Applications {\bf15}:59 no.3, 2024.}




{\small\bibitem{Hong-Lu-Zhang}

Q. Hong, G. Lu and  L. Zhang, {\it $\L^p$ Boundedness of rough bi-parameter Fourier integral operators}, Forum Math {\bf 30 }: no.1, 87-107, 2018.}






{\small\bibitem{Ruzhansky-Sugimoto}
M. Ruzhansky and M. Sugimoto, { \it A local-to-global boundedness argument and Fourier integral operators}, Journal of Mathematical Analysis and Applications {\bf 473}: no.2,  892-904, 2019.}

{\small\bibitem{Sindayigaya}	
	J. Sindayigaya, {\it On the global $L^2$-boundedness of Fourier integral operators with rough amplitude and phase functions}, Forum Mathematicum {\bf 35}: no. 3, 2023.
}

{\small\bibitem{Wu-Yang}
	
	G. N.  Wu  and J.Yang, {\it On $ L^p$ boundedness of rough Fourier integral operators}, Forum Mathematicum {\bf 36}: no.6, 2024.}




	













\end{thebibliography}
\end{document}